\documentclass[11pt]{article}
\usepackage[dvips]{color}
\usepackage{mathrsfs}
\usepackage{graphicx}
\usepackage{amsfonts}
\usepackage{amssymb}
 \usepackage{epstopdf}
\definecolor{dgreen}{rgb}{0,.8,.3}
\definecolor{blue}{rgb}{.2,.3,.7}
\definecolor{red}{rgb}{1.0,0.2,0.2}
\usepackage{color}
\input amssym.def

\input epsf.tex

\begin{document}

\renewcommand{\Box}{\rule{2.2mm}{2.2mm}}
\newcommand{\BOX}{\hfill \Box}

\newtheorem{eg}{Example}[section]
\newtheorem{thm}{Theorem}[section]
\newtheorem{lemma}{Lemma}[section]
\newtheorem{example}{Example}[section]
\newtheorem{remark}{Remark}[section]
\newtheorem{proposition}{Proposition}[section]
\newtheorem{corollary}{Corollary}[section]
\newtheorem{defn}{Definition}[section]
\newtheorem{alg}{Algorithm}[section]
\newtheorem{ass}{Assumption}[section]
\newenvironment{case}
    {\left\{\def\arraystretch{1.2}\hskip-\arraycolsep \array{l@{\quad}l}}
    {\endarray\hskip-\arraycolsep\right.}

\def\argmin{\mathop{\rm argmin}}

\makeatletter
\renewcommand{\theequation}{\thesection.\arabic{equation}}
\@addtoreset{equation}{section} \makeatother

\title{Relaxed constant positive linear dependence constraint qualification and its application to bilevel programs}
\author{ Mengwei Xu\thanks{\baselineskip 9pt School of Mathematics, Tianjin University, Tianjin, 300072, China.
E-mail: xumengw@hotmail.com. The
research of this author was supported by the National Natural
Science Foundation of China under Projects  No. 11601376} \ and \   Jane J. Ye\thanks{\baselineskip 9pt Corresponding author. Department of Mathematics
and Statistics, University of Victoria, Victoria, B.C., Canada V8W 2Y2. E-mail: janeye@uvic.ca.
The research of this author was partially supported by NSERC.}
}
\date{}
\maketitle

\baselineskip 18pt

{\bf Abstract.} Relaxed constant positive linear dependence constraint qualification (RCPLD) for a system of smooth equalities and inequalities  is a   constraint qualification that is weaker than the usual constraint qualifications such as Mangasarian Fromovitz constraint qualification and the linear constraint qualification. Moreover RCPLD is known to induce an error bound property. In this paper we extend  RCPLD to a very general  feasibility system  which may include  Lipschitz continuous  inequality constraints, complementarity constraints and {abstract} constraints. We show that   {this RCPLD  for the general system is  a constraint qualification for the optimality  condition in terms of limiting subdifferential and limiting normal cone} and it is a sufficient condition for the error bound property under the strict complementarity condition for the complementarity system and Clarke regularity conditions for the inequality constraints and the abstract constraint set. Moreover we introduce and study some sufficient conditions for RCPLD including the relaxed constant rank constraint qualification  (RCRCQ).  Finally we apply our results to  the bilevel program.



\newpage

\baselineskip 18pt
\parskip 2pt

\section{Introduction}
A constraint qualification is a condition imposed on the constraint region of a mathematical program under which the
 Karush-Kuhn-Tucker (KKT) condition holds at any local optimal solution. Other than guaranteeing  KKT condition  holding at all local optimal solutions, some constraint qualifications also lead to  existence of error bounds to the feasible region and hence play a key role in
 convergence analysis  of certain computational methods. Hence studying constraint qualifications is essential in both theoretical and numerical points of view.

For smooth nonlinear programs with equality and inequality constraints, the classical constraint qualifications are the linear independence constraint qualification (LICQ),
Mangasarian-Fromovitz constraint qualification (MFCQ), and  the linear constraint qualification (LCQ), i.e., all functions in the equality and inequality constraints are affine. These three classical constraint qualifications may be too restrictive for many problems.
Janin \cite{crcq} relaxed  LICQ and proposed the constant rank constraint qualification (CRCQ),
which neither implies nor is implied by MFCQ. {The concept of CRCQ was weakened to  the relaxed constant rank constraint qualification (RCRCQ) which is shown to be a constraint qualification by Minchenko and Stakhovshi \cite{M-S}.}
Qi and Wei \cite{CPLD} introduced the concept of the constant positive linear dependence (CPLD) condition which is weaker than both  CRCQ and  MFCQ. CPLD was  shown to be a constraint qualification by Andreani, Mart\'{i}nez and Schuverdt  in \cite{ams}.

In \cite{RCPLD}, Andreani,   Haeser, Schuverdt and Silva introduced the  following relaxed version of  CPLD  for a system of smooth equality and inequality constraints:
$$ g_i(x)\leq 0,\ i=1,\cdots,n,\ \ 
 h_i(x)=0,\ i=1,\cdots,m,$$
 where $g_i, h_i:\mathbb{R}^d\rightarrow \mathbb{R}$ are smooth at  $x^*$,  a feasible solution.   $x^*$ is said to satisfy  the relaxed constant positive linear dependence constraint qualification $(\rm RCPLD)$ if there exists $\mathbb{U}(x^*)$, a neighborhood  of $x^*$ such that
 \begin{itemize}
\item[{\rm (i)}]   $\{\nabla h_i(x)\}_{i=1}^m$ has the same rank for every $x\in \mathbb{U}(x^*)$.
\item[{\rm (ii)}]   Let $J\subseteq \{1,\cdots,m\}$ be such that $\{\nabla h_i(x^*)\}_{i\in J}$ is a basis for  span$\{\nabla h_i(x^*)\}_{i=1}^m$.  For every $I\subseteq I_g^*:=\{i: g_i(x^*)=0\}$, if $\{\nabla g_i(x^*)\}_{i\in I}\cup 
\{\nabla h_i(x^*)\}_{i\in J}$ is positive linearly dependent, i.e.,
there exist scalars $\lambda_i\geq 0$,  $i\in I$, $\mu_i$, $i\in J$ not all zero  such that
\begin{eqnarray*}
&& 0=  \sum_{i\in I} \lambda_i\nabla g_i(x^*)+\sum_{i\in J} \mu_i\nabla h_i(x^*),
\end{eqnarray*}
then $\{\nabla g_i(x)\}_{i\in I} \cup \{\nabla h_i( x)\}_{i\in J}$ is linearly dependent for every $x \in \mathbb{U}(x^*)$.
\end{itemize}
It is easy to see that in the case where either  LCQ or  MFCQ holds,  RCPLD also holds. Hence  RCPLD is  weaker than LCQ and MFCQ.
In \cite{RCPLD}, the authors not only showed that  RCPLD is a constraint qualification but also proved that if all functions $g_i,h_i$ have second-order derivatives at all points near the point $x^*$, then RCPLD is a sufficient condition for  the error bound property: $\exists \alpha>0, \mathbb{U}(x^*)$, a neighborhood of $x^*$ such that
$$
d_{{\cal F}}(x) \leq \alpha (\|g_+(x)\| + \|h(x)\| ), \quad \forall x \in \mathbb{U}(x^*),
$$
where ${{\cal F}}:=\{x|g(x)\leq 0, h(x)=0\}$, $g(x)=(g_1(x),\dots, g_n(x))$, $h(x)=(h_1(x),\dots, h_m(x))$, $d_{{\cal F}}(x)$ is the distance from $x$ to set ${{\cal F}}$,  $\|\cdot\|$ denotes any norm  and $g_+(x):=\max\{0, g(x)\}$, where the maximum is taken component-wise. Moreover as an open question, in \cite{RCPLD},  a question was asked on whether or not it was possible to prove the error bound property without imposing the second-order differentiability of all functions. In Guo, Zhang and Lin \cite{gzl}, it was shown that the error bound property holds under  RCPLD without imposing the second order differentiability of all functions. Other than using it as a constraint qualification to ensure the KKT condition holds,  RCPLD
 is also used in  the convergence analysis of the augmented Lagrangian method to obtain a KKT point (see e.g., \cite{RCPLD,cglj,isu,wy15}).
Recently, Guo and Ye \cite[Corollary 3]{GuoYe} extended  RCPLD to the case where there is an extra abstract constraint set and showed that it is still a constraint qualification.  In \cite[Definition 4.3]{gly1}, a version of  RCPLD called MPEC RCPLD was introduced  for the mathematical programs with
equilibrium constraints (MPEC) and was shown in \cite[Corollary 4.1]{gl} that it is a constraint qualification for M-stationary conditions. Moreover, it was shown in \cite[Theorem 5.1]{gzl} that the RCPLD for MPECs ensures the error bound property under the assumption of strict complementarity.

In this paper, we extend  RCPLD to the following very general feasibility system:
\begin{equation}
\begin{array}{l}
g_i(x)\leq 0,\ i=1,\cdots,n,\\
 h_i(x)=0,\ i=1,\cdots,m,\\
(G(x) ,H(x)) \in \Omega^p,\\
x\in C,
 \end{array} \label{feasibilitys}
\end{equation}
where $\Omega^p:=\{ (y,z)\in \mathbb{R}^p \times \mathbb{R}^p| 0\leq y\perp z\geq 0\}$ is the $p$th dimensional  complementarity set,
$C:=C_1\times C_2\times \cdots \times C_l$ with 
$C_i\subseteq \mathbb{R}^{q_i}$  closed, $i=1,\cdots,l$, $q_1+\cdots+q_l=d$,
 the functions $ g_i:\mathbb{R}^d\to \mathbb{R},  i=1,\cdots,n$, are locally Lipschitz continuous and $h_i:\mathbb{R}^d\to \mathbb{R},  i=1,\cdots,m
$, $G, H:\mathbb{R}^d\rightarrow \mathbb{R}^p$, are continuously differentiable at $x^*$, a feasible solution. 

Denote   the feasible region of system $(\ref{feasibilitys})$ by ${\cal F}$.
For a feasible point $x^*\in {\cal F}$, we define the following index sets:
 \begin{eqnarray*}
&& I_g^*:=\{i=1,\cdots,n: g_i(x^*)=0\},\\
 && \mathcal{I}^*:=\{i=1,\cdots,p: 0=G_i(x^*)< H_i(x^*) \},\\
  && \mathcal{J}^*:=\{i=1,\cdots,p: 0=G_i(x^*)=H_i(x^*) \},\\
   && \mathcal{K}^*:=\{i=1,\cdots,p: G_i(x^*)>H_i(x^*) =0\}.
\end{eqnarray*}

\begin{defn}[RCPLD for the  nonsmooth system (\ref{feasibilitys})]\label{nrcpld} We say that  the relaxed constant positive linear dependence constraint qualification $(\rm RCPLD)$ holds  at $x^*\in {\cal F}$ for system (\ref{feasibilitys}) if
the following conditions hold:
\begin{itemize}
\item[{\rm (i)}]
The vectors $\{\nabla h_i(x)\}_{i=1}^m\cup \{\nabla G_i(x)\}_{i\in \mathcal{I}^*}\cup \{\nabla H_i(x)\}_{i\in \mathcal{K}^*}$
 have the same rank 
 for all $x$ in a  neighbourhood of $x^*$.
\item[{\rm (ii)}]   Let  $\mathcal{I}_1\subseteq \{1,\cdots,m\}$, $\mathcal{I}_2\subseteq \mathcal{I}^*$,  $\mathcal{I}_3\subseteq \mathcal{K}^*$ be such that  the set of vectors $\{\nabla h_i(x^*)\}_{i\in \mathcal{I}_1}\cup \{\nabla G_i(x^*)\}_{i\in \mathcal{I}_2}\cup \{\nabla H_i(x^*)\}_{i\in \mathcal{I}_3}$ 
 is a basis for  
 $$ \mbox{ span }\large \{ \{\nabla h_i(x^*)\}_{i=1}^m\cup \{\nabla G_i(x^*)\}_{i\in \mathcal{I}^*}\cup \{\nabla H_i(x^*)\}_{i\in \mathcal{K}^*} \large \}.$$
  For any index sets $\mathcal{I}_4\subseteq I_g^*$, $\mathcal{I}_5,\mathcal{I}_6\subseteq \mathcal{J}^*$, 
if  there exists a nonzero vector $(\lambda^g,\lambda^h, \lambda^G,\lambda^H, \eta^*) \in \mathbb{R}^{n}\times\mathbb{R}^{m}\times \mathbb{R}^{p}\times\mathbb{R}^p\times\mathbb{R}^d$ satisfying $\lambda_i^g\geq 0$ for $i\in \mathcal{I}_4$, 
and $\mbox{either } \lambda_i^G> 0, \lambda_i^H>  0 \mbox{ or }\lambda_i^G \lambda_i^H=0, \forall i \in \mathcal{J}^*$,  ${\eta^*=(\eta_1^*,\cdots,\eta_l^*)}\in \mathcal{N}_{C}(x^*)=\mathcal{N}_{C_1}(x_1^*)\times \cdots\times \mathcal{N}_{C_l}(x_l^*)$, $v^*_i \in \partial g_i(x^*)$ for $i\in {\cal I}_4$  such that 
\begin{eqnarray}\label{nna}
&&0 = \sum_{i\in \mathcal{I}_4} \lambda_i^g v_i^* +\sum_{i\in \mathcal{I}_1} \lambda_i^h \nabla h_i(x^*)
 -\sum_{i\in\mathcal{I}_2\cup\mathcal{I}_5} \lambda_i^G \nabla G_i(x^*)
 -\sum_{i\in\mathcal{I}_3\cup\mathcal{I}_6}\lambda_i^H \nabla H_i(x^*)
 + \eta^*,\nonumber
\end{eqnarray} 
then for all $k$  sufficiently large,  the set of vectors 
\begin{equation}\{v_i^k\}_{i\in \mathcal{I}_4}\cup \{\nabla h_i(x^k)\}_{i\in \mathcal{I}_1}\cup \{\nabla G_i(x^k)\}_{i\in \mathcal{I}_2\cup\mathcal{I}_5}\cup\{\nabla H_i(x^k)\}_{i\in \mathcal{I}_3\cup\mathcal{I}_6}\cup 
{ \{ \nu_i^k\}_{i\in L} }, \label{normalvector}
\end{equation} 
where $v_i^k \in \partial g_i(x^k)$, $\nu_i^k:=\{0\}^{s_i}\times \{\eta_i^k\}\times \{0\}^{t_i} \mbox{ with } \eta_i^k \in \mathcal{N}_{C_i}(x_i^k)$,  $s_i:=q_1+\cdots +q_{i-1}$, $t_i:=q_{i+1}+\cdots + q_l$, $L:= \{1,\dots, l:\eta_i^k\neq 0\}$  and $x^k\neq x^*$,
is linearly dependent for all  sequences $\{x^k\}, \{v^k\}, \{\eta^k\}$ satisfying $x^k\rightarrow x^*$,  $v_i^k\rightarrow v_i^*$,  $\eta^k:=(\eta_1^k,\dots, \eta_l^k)=\sum_{i\in L} \nu_i^k \rightarrow \eta^*$
 as $k\rightarrow \infty$.
\end{itemize}
\end{defn}
\begin{remark}\label{remark} { In Definition \ref{nrcpld} (ii), we use  sequences instead of neighborhoods. Although we could also use a neighborhood in the definition equivalently, it is more convenient to use the sequential form since if the point $x^*$ is an isolated non-differentiable point, i.e., there exists a neighborhood around $x^*$ where $g$ is differentiable, then $v_i^k
$ can be taken as the gradient $\nabla g_i(x^k)$.   Since a Lipschitz continuous function is differentiable almost everywhere and so such points are abundant. }

 In the case where the rank of $\{\nabla h_i(x^*)\}_{i=1}^m\cup \{\nabla G_i(x^*)\}_{i\in \mathcal{I}^*}\cup \{\nabla H_i(x^*)\}_{i\in \mathcal{K}^*}$  is equal to $d$, it is easy to see that RCPLD holds automatically.
What is more, in Theorem \ref{err2} we will show that the error bound property  holds in this case. 

Note that Definition \ref{nrcpld} is weaker than the one defined  {in} Guo and Ye \cite[Corollary 3]{GuoYe} for the system containing only smooth equality and inequality constraints and one  abstract constraint, in which the stronger condition $\{\nabla g_i(x)\}_{i\in \mathcal{I}_4} \cup \{\nabla h_i( x)\}_{i\in \mathcal{I}_1}$ is linearly dependent for every $x \in \mathbb{U}(x^*)$ required.

If the set of vectors in (\ref{normalvector}) is replaced by the following set of vectors
$$\{v_i^k\}_{i\in \mathcal{I}_4}\cup \{\nabla h_i(x^k)\}_{i\in \mathcal{I}_1}\cup \{\nabla G_i(x^k)\}_{i\in \mathcal{I}_2\cup\mathcal{I}_5}\cup\{\nabla H_i(x^k)\}_{i\in \mathcal{I}_3\cup\mathcal{I}_6}\cup 
{ \{ \eta^k \} },$$
then since $\eta^k=\sum_{i\in L} \nu_i^k$,
 the resulting condition would be stronger than the RCPLD defined in Definition \ref{nrcpld}. We illustrate this by Example \ref{nonsmoothex}.


\end{remark}

%

 In this paper we show that   RCPLD for the nonsmooth feasibility system is a constraint qualification for any optimization problem with a Lipschitz objective function and the constraints described as in the feasibility system (\ref{feasibilitys}). Moreover with some extra  conditions, we will show that RCPLD is a
 sufficient condition for the following error bound property: $\exists \alpha>0, \mathbb{U}(x^*)$, a neighborhood of $x^*$ such that
\begin{equation}
d_{{\cal F}}(x) \leq \alpha (\|g_+(x)\| + \|h(x)\| +\sum_{i=1}^pd_{\Omega}(G_i(x), H_i(x))), \quad \forall x \in \mathbb{U}(x^*)\cap C, \label{errorbp}
\end{equation}
 where $\Omega:=\Omega^1=\{ (y,z)\in \mathbb{R}^2| \ 0\leq y\perp z\geq 0\}$.
 
One of the motivations to extend the concept of RCPLD to the nonsmooth system (\ref{feasibilitys})  is to study the constraint qualification and optimality condition for the following bilevel program:
\begin{eqnarray*}
({\rm BP})~~~~~~\min && F(x,y)\nonumber\\
{\rm s.t.} && y\in S(x), \ 
 G(x,y)\leq 0,\ H(x,y)= 0,
\end{eqnarray*}
where  $S(x)$ denotes the solution set of the lower level program
\begin{eqnarray*}
({\rm P}_x)~~~~~~~~
\min_{y\in Y(x)} \  f(x,y),
\end{eqnarray*}
and $Y(x):=\{y\in \mathbb{R}^s: g(x,y)\leq 0, h(x,y)=0\}$, 
 $F:\mathbb{R}^d\times \mathbb{R}^s \rightarrow \mathbb{R}$ is locally Lipschitz continuous, $G:\mathbb{R}^d\times \mathbb{R}^s \rightarrow \mathbb{R}^p$ and  $H:\mathbb{R}^d\times \mathbb{R}^s \rightarrow \mathbb{R}^q$  are continuously differentiable, 
 $f:\mathbb{R}^d\times \mathbb{R}^s \rightarrow \mathbb{R}$, $g:\mathbb{R}^d\times \mathbb{R}^s \rightarrow \mathbb{R}^m, h:\mathbb{R}^d\times \mathbb{R}^s \rightarrow 
 \mathbb{R}^n$ are continuously differentiable and twice continuously differentiable with respect to variable $y$.

Bilevel programs naturally fall in the domain of global optimization since in the lower level problem, the global optimal solution is always required in either  optimality conditions or  numerical algorithms. A popular approach to deal with (BP) is to replace the set of global solutions $S(x)$ by the KKT optimality conditions of the lower level problem. This reformulation is based on the fact that if the lower level problem $(P_x)$ is a convex program for each fixed $x$ and certain constraint qualification holds, then $y\in S(x)$ if and only if  its KKT optimality condition holds. But due to the introduction of multipliers for the lower level problem as extra variables,  the resulting reformulation may not  be equivalent to the original bilevel program even in the case where $(P_x)$ is convex but the multipliers are not unique. For the discussion of this issue and the recent new results, the reader is referred to recent paper \cite{ye2019} and the reference within for more discussions.
If the lower level problem $(P_x)$ is not a convex program for each fixed $x$, the KKT condition for the lower level problem is usually only a necessary but not sufficient condition for optimality and moreover, as pointed out by Mirrlees in \cite{Mirrlees99}, such a reformulation by the KKT condition may miss out the true optimal solution of the original bilevel program.

Instead of using  $y\in S(x)$ as a constraint in (BP), Outrata \cite{Outrata} proposed  to replace it with $f(x,y)- V(x)=0, y\in Y(x)$ where $V(x) :=\inf_{y\in Y(x)} f(x,y)$ is the value function of the lower level program for a numerical consideration. This so-called  value function approach was further used in Ye and Zhu \cite{yz}  and later in other papers such as \cite{Dam-Dut-Mor,Dam-Zem,m1,m,mb,MorNamPhan} to derive various necessary optimality conditions under the partial calmness condition. The partial calmness condition for the value function reformulation of the bilevel program, however, is a very strong assumption. To derive a necessary optimality condition under weaker assumptions, Ye and Zhu  \cite{yz2}  proposed the following combined program where both the value function constraint and the KKT condition of the lower level program are used as constraints:
\begin{eqnarray*}
({\rm CP})~~~~~~\min_{x,y,u,v} && F(x,y)\nonumber\\
{\rm s.t.} && f(x,y)-V(x) \leq 0,\ G(x,y)\leq 0, H(x,y)= 0,\nonumber\\
&&\nabla_y f(x,y)+   \nabla_y g(x,y)^Tu+ \nabla_y h(x,y)^Tv=0,\\
&& (-g(x,y),u)\in \Omega^m.
\end{eqnarray*}
 As discussed in  \cite{yz2}, such a reformulation can avoid some difficulties caused by using the value function or the classical KKT approach alone. In \cite{y11}, necessary optimality conditions  in the form of  Mordukhovich (M-) 
stationary condition for (CP) is studied under the partial calmness/weak calmness condition. These conditions, however, may not be easy to verify.
 
Note that the inclusion of the value function makes the problem nonsmooth since the value function is usually nonsmooth but under some reasonable assumptions on the lower level problem, the value function is Lipschitz continuous. Hence the feasible set of (CP) is a special case of the general feasibility system (\ref{feasibilitys}). However as it was shown in Ye and Zhu \cite[Proposition 1.3]{yz2}, the no nonzero abnormal multiplier constraint qualification (NNAMCQ) as defined in Definition \ref{nna*} never holds for (CP).  Being able to deal with nonsmooth inequality constraints in  RCPLD would allow us to present verifiable constraint qualifications and exact penalty for the reformulation of  the bilevel program in which the value function is used, such as (CP).

The rest of the paper is organized as follows.  In Section 2, we present basic definitions
as well as some preliminaries which will be used in this paper. In Section 3, we show that  RCPLD  is a constraint qualification and a sufficient condition for error bound properties under certain regularity conditions.  We introduce some sufficient conditions for RCPLD, which are easier to verify in Section 4. In Section 5, we apply  RCPLD to  the bilevel program.

We adopt the following standard notation in this paper. For any two vectors $a$ and $b $, we denote by either $\langle a, b \rangle$ or $a^T b$  their inner product.  Given a  differentiable function $G: \mathbb{R}^d\rightarrow \mathbb{R}$, we denote its
gradient vector by $\nabla G(x)\in \mathbb{R}^d$.  For a differentiable  mapping $\Phi:\mathbb{R}^d\to \mathbb{R}^n$ with $n\geq 2$ and a vector $x\in \mathbb{R}^d$, we denote by
  $\nabla \Phi(x)\in \mathbb{R}^{n\times d}$ the Jacobian matrix of $\Phi$ at $x$. For a set $C\subseteq \mathbb{R}^d$, we denote by int $C$, co $C$, $\bar C$, bd $C$ and $d_C(x)$ the interior, the convex hull, the closure, the boundary of $C$ and the distance from $x$ to $C$, respectively. 
We denote by $|{\cal I}|$ the cardinality of index set  ${\cal I}\subseteq \{1,\dots, n\}$.  For any vector $v\in \mathbb{R}^n$ and a given index set ${\cal I}\subseteq \{1,\dots, n\}$, we use $v_{\cal I}$ to denote the vector in $ \mathbb{R}^{|{\cal I}|}$ with components $v_i, i\in {\cal I}$.  For a matrix $A\in \mathbb{R}^{n\times m}$, $A^T$ denotes its transpose, $r(A)$ denotes its rank.  We denote by ${\mathbb{B}}_\delta(\bar x)$  the open  ball center at $\bar x$  with radius  $\delta>0$ and by ${\mathbb{B}}$ the open unit ball center at the origin. Unless otherwise specified we denote by $\|\cdot\|$ any  norm in the finite dimensional space.

\section{Background and preliminaries}
In this section, we present some background materials on variational analysis which will be used throughout the paper.
Detailed discussions on these
subjects
 can be found in \cite{c,clsw,m1,var}.

 Given a set  $C\subseteq \mathbb{R}^d$, its Fr\'{e}chet/regular normal cone at $z\in C$ is  defined by
$$
\widehat{\mathcal{N}}_{C}(z):= \{v \in \mathbb{R}^d:  \langle v,z'-z\rangle \le
 o\big(\|z-z'\|\big)\;\mbox{ for all }\; z'\in C\}
 $$
and the limiting/Mordukhovich/basic normal cone  to $C$ at point $ z$ is defined by
$$
\mathcal{N}_{C}(z)=\{\lim_{k\rightarrow \infty} v_k: v_k \in \widehat{\mathcal{N}}_C(z_k) ,\ \ z_k \in C,   z_k\rightarrow z\}.
$$
A set $C\subseteq \mathbb{R}^d$ is Clarke regular at $z$ if it is locally closed at $z$ and $\mathcal{N}_{C}(z)=\widehat{\mathcal{N}}_{C}(z)$ \cite[Definition 6.4]{var}. Note that for $C:=C_1\times C_2\times \cdots \times C_l$ with 
$C_i$ closed, $i=1,\cdots,l$, by \cite[Proposition 6.41]{var}, $C$ is regular at $z$ if and only if each $C_i$ is regular at $z_i$.

The exact expressions for the limiting normal cone of the complementarity set present below will be useful. In this paper we denote by $\Omega^p$ the  {complementarity set} and when $p=1$,  $\Omega:=\Omega^1=\{ (y,z)\in \mathbb{R}^2| \ 0\leq y\perp z\geq 0\}$.
\begin{proposition}(See e.g.  \cite[Proposition 3.7]{y00}) \label{normalC}
For any $(a,b)$ lying in the complementarity set $ \Omega^p$, the limiting normal cone to $\Omega^p$ at $(a,b)$ is
\begin{eqnarray*}
\mathcal{N}_{\Omega^p} (a,b)=\left\{-(\alpha,\beta) \in \mathbb{R}^p\times \mathbb{R}^p :
\begin{array}{ll}
\alpha_i \in \mathbb{R},\ \beta_i=0  & {\rm if}\ a_i=0<b_i\\
\alpha_i=0,\ \beta_i \in \mathbb{R}  & {\rm if}\ a_i>0=b_i\\
{\rm either}\ \alpha_i>0,\beta_i>0\ {\rm or}\ \alpha_i \beta_i=0\  & {\rm if}\ a_i=b_i=0
 \end{array}\right\}.
 \end{eqnarray*}
\end{proposition}



  Let $f: \mathbb{R}^d \rightarrow \overline \mathbb{R}$ be a lower semicontinuous function and finite at $ x\in \mathbb{R}^d$.
 We define the Fr\'{e}chet/regular  subdifferential
(\cite[Definition 8.3]{var})
 of $f $ at $  x$
as
\begin{eqnarray*}
\widehat{\partial} f( x)
:=
\{\zeta \in \mathbb{R}^d:  f(x') - f( x) -\langle \zeta, x'-x \rangle  \geq o\|x'- x\| \},
\end{eqnarray*}
and the limiting/Mordukhovich/basic  subdifferential of $f$ at ${x}$
as
\begin{eqnarray*}
\partial f({x}):=
\{\lim_{k\rightarrow \infty} \xi_k: \xi_k \in \widehat{\partial}f(x_k), x_k\rightarrow  x,f(x_k)\rightarrow f(x)\}.
\end{eqnarray*}
Let $f: \mathbb{R}^d \rightarrow  \mathbb{R}$ be Lipschitz continuous at $ x$. We say that $f$ is subdifferentially/Clarke regular at $ x$
provided  that $\partial f({x})=\widehat{\partial} f(  x)$ \cite[Corollary 8.11]{var}.



 The following proposition collects some useful properties and calculus rules of the limiting subdifferential.

\begin{proposition} \label{sumrule2}
\begin{itemize}
\item[\rm (i)] {\rm \cite[Exercise 10.10]{var} and \cite[Theorem 2.33]{m1}}  Let $f, g:\mathbb{R}^n \to [-\infty,\infty]$ be proper lower semi-continuous around $x^* \in \mathbb{R}^n$ and finite at $x^*$, and $\alpha,\beta$ be nonnegative scalars.  Assume that at least one of them is Lipschitz continuous around $x^*$. Then
$$ \partial (\alpha f+\beta g)(x^*) \subseteq \alpha \partial f(x^*)+\beta \partial  g(x^*).$$
\item[\rm (ii)]{\rm \cite[Theorem 2.5, Remark (2)]{jourani1993} and \cite[Theorem 3.41]{m1}}Let $g: \mathbb{R}^n\to\mathbb{R}^m$ be Lipschitz continuous around $x^*$ and $f:\mathbb{R}^m\to\mathbb{R}$ be Lipschitz continuous near $g(x^*)$. Then the composite function $f\circ g$ is Lipschitz continuous around $x^*$ and
$$\partial  (f\circ g)(x^*) \subseteq \bigcup_{\xi \in \partial  f(g(x^*))} \partial  \langle \xi, g\rangle (x^*).$$
\item[{\rm (iii)}]{\rm \cite[Theorem 3.46 and Proposition 1.113]{m1}} Let $\varphi_i : \mathbb{R}^n\to \mathbb{R}\ (i=1,\dots, n)$ be  Lipschitz continuous  around $x^*$  and
$\varphi_{\max}({x}) := \max \{\varphi_i({x}) : i = 1, \ldots, n\}, $ and 
$\varphi_{\min}({x}) := \min \{\varphi_i({x}) : i = 1, \ldots, n\}.$ Then $\varphi_{\max}(x)$ and $\varphi_{\min}(x)$  are Lipschitz continuous
around $x^*\in  \mathbb{R}^n$ and 
\begin{eqnarray*}
\partial  \varphi_{\max}(x^*) &\subseteq & co\{\partial
\varphi_i(x^*): i \in I_+(x^*)\},\\
\partial  \varphi_{\min}(x^*) &\subseteq & \{\partial
\varphi_i(x^*): i \in I_-(x^*)\},
\end{eqnarray*}  where $ I_+(x^*):= \{i :
\varphi_i(x^*) = \varphi_{\max}(x^*)\}$ and $ I_-(x^*):= \{i :
\varphi_i(x^*) = \varphi_{\min}(x^*)\}$,
and the first inclusion holds as an equation if each $\varphi_i$ is subdifferentially regular at $x^*$.
 \end{itemize}
\end{proposition}
Taking into account that all norms in a finite dimensional space are equivalent, the following formula for distance function to the complementarity set $\Omega$ can be used in the error bound property (\ref{errorbp}).
\begin{proposition} (see e.g. \cite{KS}) When the norm is chosen to be the $l_1$-norm in the distance function $d_\Omega$, for any $(a,b) \in \mathbb{R}^2$,
 $$ d_{\Omega}(a,b)=\max\{-a, -b,-(a+b),\min \{ a, b\}\}.$$
 When the norm is chosen to be the $l_\infty$-norm in the distance function $d_\Omega$, for any $(a,b) \in \mathbb{R}^2$,
 $$ d_{\Omega}(a,b)=|\min \{ a, b\}|.$$
\end{proposition}

In Theorem \ref{OC1}, we need to calculate $\partial\phi_0(x)$, where  
\begin{eqnarray*}
\phi_0(x):= \sum_{i=1}^n \max   \{0,  g_i(x)\}+\sum_{i=1}^m |h_i(x)|+\sum_{i=1}^p |\min \{G_i(x), H_i(x)\}| .\end{eqnarray*}
 In order to calculate $\partial\phi_0(x)$, we  define the following index sets for each $x$:
\begin{eqnarray}
&&A(x):=\{i=1,\cdots,n:g_i(x)\geq  0\},\nonumber\\
 && \mathcal{I}(x):=\{i=1,\cdots,p: G_i(x)<H_i(x)\},\nonumber\\
  && \mathcal{J}(x):=\{i=1,\cdots,p: G_i(x)= H_i(x) \},\label{index}\\
   && \mathcal{K}(x):=\{i=1,\cdots,p: G_i(x)>H_i(x) \}.\nonumber
\end{eqnarray}
From the calculus  rules in  Proposition \ref{sumrule2}, 
we have the following estimate for the limiting subdifferential of $\phi_0$.

\begin{lemma}\label{subphi} 
Assume that  the functions $ g_i:\mathbb{R}^d\to \mathbb{R},  i=1,\cdots,n$, are locally Lipschitz continuous and $h_i:\mathbb{R}^d\to \mathbb{R},  i=1,\cdots,m
$, $G, H:\mathbb{R}^d\rightarrow \mathbb{R}^p$, are continuously differentiable around $x^*$.
$\phi_0(x)$ is a Lipschitz continuous function and for any $x^*$, there exist $\lambda_{i}^g\geq 0, i\in A(x^*)$,
$\lambda_{i}^h,i=1,\cdots,m$ and $\lambda_{i}^G, \lambda_{i}^H$, $i=1,\cdots,p$ satisfying \begin{equation}
\left \{\begin{array}{ll}
\lambda_i^H=0 &\mbox{ if } i \in {\mathcal{I}(x^*)}\\
\lambda_i^G=0 &\mbox{ if } i \in {\mathcal{K}(x^*)}\\
\mbox{either } \lambda_i^G>0, \lambda_i^H>0 \mbox{ or }\lambda_i^G\lambda_i^H=0 &  \mbox{ if } i \in {\mathcal{J}(x^*)}
\end{array}
\right .  \label{multipliers}
\end{equation}
 such that
\begin{eqnarray*}
\partial \phi_0(x^*)\subseteq
\sum_{i\in A(x^*)}\lambda_i^g \partial g_i(x^*)+  \sum_{i=1}^m\lambda_i^h \nabla h_i(x^*)
-\sum_{i=1}^p\lambda_i^G \nabla G_i(x^*)-\sum_{i=1}^p\lambda_i^H \nabla H_i(x^*).
\end{eqnarray*}
\end{lemma}
{\bf Proof.} For $i=1,\cdots,p$, let $F_i(x):=|f_i(x)|$ and $f_i(x):=\min \{G_i(x), H_i(x)\}$. From the chain rule in Proposition \ref{sumrule2}(ii), we have
$\partial F_i(x^*)\subseteq  \{\partial(\mu f_i)(x^*):\mu\in \partial |\cdot|(f_i(x^*))\}.$
We divide the analysis into three cases.
\begin{itemize}
\item[Case 1:]  $i\in \mathcal{I}(x^*)$. In this case  we have $\ G_i(x^*)<H_i(x^*)$ and hence  $f_i(x^*)=G_i(x^*)$.
By the chain rule we have $\partial F_i(x)=\{\nabla G_i(x^*)\}$ if $G_i(x^*)>0$ and $\partial F_i(x)=\{-\nabla G_i(x^*)\}$ if $G_i(x^*)<0$. If $G_i(x^*)=0$, 
then $\partial |\cdot|(f_i(x^*))=[-1,1]$ and hence $\partial F_i(x)\subseteq \{\mu\nabla G_i(x^*):\mu\in[-1,1]\}$.

\item[Case 2:] $i\in \mathcal{K}(x^*)$. Similarly as in Case 1, we can show  $\partial F_i(x)\subseteq \{\mu\nabla H_i(x^*):\mu\in[-1,1]\}$.

\item[Case 3:] $i\in \mathcal{J}(x^*)$. In this case, $f_i(x^*)=G_i(x^*)=H_i(x^*)$. If $G_i(x^*)=H_i(x^*)<0$, then  $\partial |\cdot|(f_i(x^*))=\{-1\}$ and
$\partial F_i(x^*)\subseteq \partial (-f_i)(x^*)$. Since $-f_i(x)=\max \{-G_i(x), -H_i(x)\}$, by Proposition \ref{sumrule2}(iii), $\partial F_i(x)= co\{-\nabla G_i(x^*),-\nabla H_i(x^*)\}$. If $G_i(x^*)=H_i(x^*)>0$, then  $\partial |\cdot|(f_i(x^*))=\{1\}$ and $\partial F_i(x^*)\subseteq \partial f_i(x^*)$. It follows by Proposition \ref{sumrule2}(iii) that
$\partial F_i(x)\subseteq \{\nabla G_i(x^*),\nabla H_i(x^*)\}$. If $G_i(x^*)=H_i(x^*)=0$, then  $\partial |\cdot|(f_i(x^*))=[-1,1]$, $\partial F_i(x)\subseteq  \{\partial(\mu f_i)(x^*):\mu\in [-1,1]\}$. If $\mu>0$, then $\partial(\mu f_i)(x^*) \subseteq \mu \partial f_i
(x^*) \subseteq \mu\{\nabla G_i(x^*),\nabla H_i(x^*)\} $. If $\mu<0$, then 
 $\partial(\mu f_i)(x^*)=-\mu \partial (-f_i)(x^*)  =-\mu co\{-\nabla G_i(x^*),-\nabla H_i(x^*)\}$. Hence 
 $$ \partial F_i(x^*) \subseteq \{\lambda v:\lambda\geq 0, v\in co\{-\nabla G_i(x^*),-\nabla H_i(x^*)\}\cup \{\nabla G_i(x^*),\nabla H_i(x^*)\}\}.$$
\end{itemize}
Summarizing the above cases, we have that for any $x^*$, 
$$\partial |\min \{G_i(x^*), H_i(x^*)\}|\subseteq  \left(\begin{array}{l}
\nabla G(x^*)\\
\nabla H(x^*)
 \end{array}\right)^T (-\lambda^G,-\lambda^H),$$ where $(\ref{multipliers})$ holds.

Since $g_i(x),i=1,\cdots,n$, $|h_i(x)|, j=1,\cdots,m$ and 
$|\min \{G_i(x), H_i(x)\}|, i=1\cdots, p$ are all Lipschitz continuous around $x^*$, by the calculus  rules in  Proposition \ref{sumrule2} and (\ref{multipliers}),
there exist parameters $\lambda_{i}^g\geq 0, i\in A(x^*)$, 
$\lambda_{i}^h,i=1,\cdots,m$ and $\lambda_{i}^G, \lambda_{i}^H$, $i=1,\cdots,p$ satisfying (\ref{multipliers}) such that 
\begin{eqnarray*}
\partial \phi_0(x^*)\subseteq
\sum_{i\in A(x^*)}\lambda_i^g \partial g_i(x^*)+  \sum_{i=1}^m\lambda_i^h \nabla h_i(x^*)
-\sum_{i=1}^p\lambda_i^G \nabla G_i(x^*)-\sum_{i=1}^p\lambda_i^H \nabla H_i(x^*).
\end{eqnarray*}  
\BOX

We will  need the following  result which is an extension of  Carath\'{e}odory's lemma.
\begin{lemma}\cite[Lemma 1]{RCPLD}\label{lem3-1}
If $v=\displaystyle \sum_{i=1}^{m+n} \alpha_i v_i$ with $v_i\in\mathbb{R}^d$ for every $i$, $\{v_i\}_{i=1}^m$ is linearly independent and $\alpha_i\neq 0$, $i=m+1,\cdots,m+n$, then there exist $I\subseteq \{m+1,\cdots,m+n\}$ and scalars $\bar{\alpha}_i$ for every $i\in \{1,\cdots,m\}\cup I$ such that\\
(i) $v=\displaystyle \sum_{\{1,\cdots,m\}\cup I} \bar{\alpha}_i v_i$ with $\alpha_i \bar{\alpha}_i>0$ for every $i\in I$,\\
(ii) $\{v_i\}_{\{1,\cdots,m\}\cup I}$ is linearly independent.
\end{lemma}

\section{RCPLD as a constraint qualification and a sufficient condition for error bounds}
We first  show that the RCPLD introduced in Definition \ref{nrcpld} is a constraint qualification for any optimization problem in the form 
$
({\rm P})~~\min_{x} \quad  f(x) \quad 
{\rm s.t.} \quad x\in \mathcal{F},
$
where $f:\mathbb{R}^d\rightarrow \mathbb{R}$ is  Lipschitz continuous around any local optimal solution and $\mathcal{F}$ is the feasible region defined by the system (\ref{feasibilitys}).

 \begin{thm}\label{OC1} Let $x^*$ be a local minimizer of the optimization problem $({\rm P})$.
Suppose that  RCPLD holds at $x^*$.
Then $x^*$ is an M-stationary point, i.e., there exist
$\lambda_{i}^g\geq 0$, $i\in I_g^*$, $\lambda_{i}^h$, $i=1,\cdots,m$, $\lambda_{i}^G$, $i=1,\cdots,p$ and $\lambda_{i}^H$, $i=1,\cdots,p$ such that 
\begin{eqnarray*}
&& 0\in\partial f(x^*)+\sum_{i\in I_g^*} \lambda_i^g \partial g_i(x^*) +\sum_{i=1}^m \lambda_i^h \nabla h_i(x^*)\\
&&~~~~-\sum_{i\in \mathcal{I}^*\cup\mathcal{J}^*}\lambda_i^G \nabla G_i(x^*)-\sum_{i\in \mathcal{K}^*\cup\mathcal{J}^*}\lambda_i^H \nabla H_i(x^*)+\mathcal{N}_{C}(x^*),\\
&&  \lambda_i^H=0 \  i\in {\cal I}^*, \lambda_i^G=0 \ i\in {\cal K}^*,  \mbox{ either } \lambda_i^G>0, \lambda_i^H>0, \mbox{ or } \lambda_i^G\lambda_i^H=0,\ i\in \mathcal{J}^*.
\end{eqnarray*}
\end{thm}
{\bf Proof.} 
Step 1:  Note that a point $x\in C$ is feasible for problem (${\rm P}$) if and only if{$$\phi_0(x):=\sum_{i=1}^n \max   \{0,  g_i(x)\}+\sum_{i=1}^m |h_i(x)|+\sum_{i=1}^ p |\min \{G_i(x), H_i(x)\}| =0.$$}
 For each $k$, we consider the following penalized problem:
\begin{eqnarray*}
({\rm P}_k)~~~\min && F^k(x):=f(x)+\frac{k}{2}\phi_0^2(x)+\frac{1}{2}\|x-x^*\|^2 \\
 {\rm s.t.} && x\in { \bar \mathbb{B}}_{\varepsilon}(x^*)\cap C,
\end{eqnarray*}
where $\|\cdot\|$ denotes the 2-norm, $\varepsilon>0$  is such that $f(x^*)\leq f(x)$ for all feasible point $x\in \bar \mathbb{B}_{\varepsilon}(x^*)$. 
Since the feasible region is compact and the objective function is continuous,  there exists an optimal solution $x^k$ of the problem $({\rm P}_k)$. Taking subsequence if necessary, we assume that $\displaystyle \lim_{k\to\infty}x^k=\bar x$ and thus $\bar x\in  \bar \mathbb{B}_{\varepsilon}(x^*)\cap C$.  Moreover by the optimality of $x^k$, we have
\begin{eqnarray}\label{th2-2-1}
f(x^k)+\frac{k}{2}\phi_0^2(x^k)+\frac{1}{2}\|x^k-x^*\|^2
= F^k(x^k) \leq F^k(x^*)=f(x^*).
\end{eqnarray}
From the boundedness of $\{f(x^k)\}$, we have that
$
\phi_0(x^k)\to 0
$
as $k\to\infty$. It follows that $\bar x$ is a feasible point of the problem $(P)$.  Condition $(\ref{th2-2-1})$ yields 
$
f(x^k)+\frac{1}{2}\|x^k-x^*\|^2 \leq f(x^*).
$
Taking limit as $k\to\infty$, we obtain 
$
f(\bar x)+\frac{1}{2}\|x^*-\bar x\|^2 \leq f(x^*),
$
which means that $\bar x=x^*$. Thus the sequence $\{ x^k\} $ converges to $x^*$. 

Step 2:  Since for sufficiently large $k$, $x^k$ is an interior point of $\bar \mathbb{B}_{\varepsilon}(x^*)$, by the necessary optimality condition in terms of limiting subdifferential and the nonsmooth calculus rule, there exist $v_0^k\in \partial f(x^k)$, $u^k\in k \phi_0(x^k) \partial\phi_0(x^k)+\eta^k$, $\eta^k
\in \mathcal{N}_{C}(x^k)$ such that
\begin{eqnarray}
0=v_0^k+u^k+ (x^k-x^*).\label{sumrule}
\end{eqnarray}

Suppose that there is a subsequence of $\{u^k\}$ such that all $u^k$ equal to zero in this subsequence.
 Since $f$ is 
Lipschitz continuous, its limiting subdifferential is compact and so the sequences  
 $\{v_0^k\}$ is compact. Taking subsequence if necessary, we assume $\displaystyle \lim_{k\rightarrow 0} u^k=0$ and
 $\displaystyle\lim_{k\to \infty}v_0^k=v_0\in \partial f(x^*) $. Taking limit as $k\to \infty$ in $(\ref{sumrule})$, we have  $0\in \partial f(x^*)$ by the outer semi-continuity of the limiting subdifferential and thus $x^*$ is a stationary point of problem (${\rm P}$) automatically. So  without loss of generality, we may assume $u^k\neq 0$ for all sufficiently large $k$.

By Lemma \ref{subphi}, since $k \phi_0(x^k)\geq 0$  there exist 
$\lambda_{i,k}^g\geq 0, i\in A(x^k)$, 
$\lambda_{i,k}^h,i=1,\cdots,m$,
 $v_i^k\in \partial g_i(x^k)$, $i\in A(x^k)$ and $\lambda_{i,k}^G, \lambda_{i,k}^H$, $i=1,\cdots,p$ such that   $\mbox{ either } \lambda_{i,k}^G> 0, \lambda_{i,k}^H>  0 \mbox{ or }\lambda_{i,k}^G \lambda_{i,k}^H=0, \ \  \forall i \in {\cal J}(x^k)$ and 
\begin{eqnarray*}
\lefteqn{k \phi_0(x^k)\partial \phi_0(x^k) 
 \subseteq  }\\
&& \sum_{i\in A(x^k)}\lambda_{i,k}^g v_i^k+\sum_{i=1}^m \lambda_{i,k}^h \nabla h_i(x^k)
 -\sum_{i\in \mathcal{I}(x^k)\cup  \mathcal{J}(x^k)  } \lambda_{i,k}^G \nabla G_i(x^k)-\sum_{i\in \mathcal{K}(x^k)\cup  \mathcal{J}(x^k) } \lambda_{i,k}^H \nabla H_i(x^k).\end{eqnarray*} 
 By the continuity of $g$, it is easy to see that  $A(x^k)\subseteq I_g^*$ for sufficiently large $k$. Similarly, by the continuity of $G$ and $H$, $\mathcal{I}(x^k) \cup \mathcal{J}(x^k) \subseteq \mathcal{I}^*\cup \mathcal{J}^*$ and $\mathcal{K}(x^k) \cup \mathcal{J}(x^k) \subseteq \mathcal{K}^*\cup \mathcal{J}^*$
for sufficiently large $k$.
Let $\lambda_{i,k}^g=0$ if $i\in I_g^* \setminus A(x^k)$,   $\lambda_{i,k}^G=0$ if $i\in  ( \mathcal{I}^*\cup \mathcal{J}^*) \setminus (\mathcal{I}(x^k) \cup \mathcal{J}(x^k)) $, and $\lambda_{i,k}^H=0$ if $i\in  ( \mathcal{K}^*\cup \mathcal{J}^*) \setminus (\mathcal{K}(x^k) \cup \mathcal{J}(x^k)) $. 
Then we have that $\lambda_{i,k}^G> 0, \lambda_{i,k}^H>  0 \mbox{ or }\lambda_{i,k}^G \lambda_{i,k}^H=0, \ \  \forall i \in {\cal J}^*$.
Since  $u^k\in k \phi_0(x^k) \partial\phi_0(x^k)+\eta^k$, it follows that 
\begin{eqnarray*}
u^k&=&\sum_{i\in I_g^*} \lambda_{i,k}^g v_i^k+\sum_{i=1}^m \lambda_{i,k}^h \nabla h_i(x^k)
-\sum_{i\in \mathcal{I}^*\cup \mathcal{J}^*} \lambda_{i,k}^G \nabla G_i(x^k)-\sum_{i\in \mathcal{K}^*\cup \mathcal{J}^*} \lambda_{i,k}^H \nabla H_i(x^k)+\eta^k.
\end{eqnarray*}
Since $\eta^k:=(\eta^k_1,\dots, \eta^k_l)\in \mathcal{N}_{C}(x^k)=\mathcal{N}_{C_1}(x^k_1)\times\cdots\times \mathcal{N}_{C_l}(x^k_l)$, we have $\nu_i^k:=\{0\}^{s_i}\times \{\eta_i^k\}\times \{0\}^{t_i}\in \mathcal{N}_{C}(x^k)$.  We denote by $\eta^k=\sum_{i\in L_k} \nu_i^k$ with $L_k:=\{i=1,\cdots,l:\eta_i^k\not =0\}$.

Let $\mathcal{I}_1\subseteq \{1,\cdots,m\}$, $\mathcal{I}_2\in \mathcal{I}^*$ and $\mathcal{I}_3\in \mathcal{K}^*$ be such that  $\{\nabla h_i(x^*)\}_{i\in \mathcal{I}_1}\cup \{\nabla G_i(x^*)\}_{i\in \mathcal{I}_2}\cup \{\nabla H_i(x^*)\}_{i\in \mathcal{I}_3}$ is a basis for $\mbox{ span } \{\{\nabla h_i(x^*)\}_{i=1}^m\cup \{\nabla G_i(x^*)\}_{i\in \mathcal{I}^*}\cup \{\nabla H_i(x^*)\}_{i\in \mathcal{K}^*}\}$.
Then by RCPLD (i), 
$\{\nabla h_i(x^k)\}_{i\in \mathcal{I}_1}\cup \{\nabla G_i(x^k)\}_{i\in \mathcal{I}_2}\cup \{\nabla H_i(x^k)\}_{i\in \mathcal{I}_3}$ is linearly independent and thus is a basis of $\mbox{ span }\{ \{\nabla h_i(x^k)\}_{i=1}^m\cup \{\nabla G_i(x^k)\}_{i\in \mathcal{I}^*}\cup \{\nabla H_i(x^k)\}_{i\in \mathcal{K}^*}\}$.
 Hence there exist 
$\tilde\lambda_{i,k}^h,i\in \mathcal{I}_1$, $\tilde{\lambda}_{k}^G, \tilde{\lambda}_{k}^H$
  such that 
\begin{eqnarray}\label{akt}
u^k &=& \sum_{i\in  I_g^*\cap {\rm supp}(\lambda_{i,k}^g)} \lambda_{i,k}^g v_i^k+\sum_{i\in \mathcal{I}_1} \tilde \lambda_{i,k}^h\nabla h_i(x^k)-\sum_{i\in \mathcal{I}_2}\tilde \lambda_{i,k}^G \nabla G_i(x^k)-\sum_{i\in \mathcal{I}_3} \tilde \lambda_{i,k}^H \nabla H_i(x^k)\nonumber \\
&& -\sum_{i\in\mathcal{J}^*\cap {\rm supp}(\tilde{\lambda}_{i,k}^G)} \tilde \lambda_{i,k}^G \nabla G_i(x^k)
-\sum_{i\in\mathcal{J}^*\cap {\rm supp}(\tilde{\lambda}_{i,k}^H)}\tilde \lambda_{i,k}^H \nabla H_i(x^k)+\sum_{i\in L_k} \nu_i^k, 
\end{eqnarray}
where ${\rm supp}(a):=\{i:a_i\neq 0\}$,  
 $\tilde \lambda_{i,k}^G>0, \tilde \lambda_{i,k}^H>0$ or $\tilde \lambda_{i,k}^G \tilde \lambda_{i,k}^H=0$ if $i\in \mathcal{J}^*$.

Since $u^k\not =0$, applying Carath\'{e}odory's lemma  in Lemma \ref{lem3-1} to  (\ref{akt}), we obtain subsets 
\begin{eqnarray*}
\mathcal{I}_4^k\subseteq I_g^*\cap {\rm supp}(\tilde{\lambda}_{i,k}^g),\qquad
\mathcal{I}_5^k\subseteq\mathcal{J}^*\cap {\rm supp}(\tilde{\lambda}_{i,k}^G),\quad \mathcal{I}_6^k\subseteq\mathcal{J}^*\cap {\rm supp}(\tilde{\lambda}_{i,k}^H), \quad L'_k\subseteq L_k,
\end{eqnarray*}
and $\{\bar{\lambda}_{k}^g, \bar{\lambda}_{k}^h, \bar{\lambda}_{k}^G,\bar{\lambda}_{k}^H, \alpha^k\}$ 
with  $\alpha_i^k> 0$, $i\in L'_k$, 
$\bar{\lambda}_{i,k}^g> 0$, $i\in \mathcal{I}_4^k$ and $\bar{\lambda}_{i,k}^G>0, \bar{\lambda}_{i,k}^H>0$ or $\bar{\lambda}_{i,k}^G \bar{\lambda}_{i,k}^H=0$ if  $i\in \mathcal{J}^*$
such that 
\begin{eqnarray}
u^k= \sum_{i\in \mathcal{I}_4^k} \bar{\lambda}_{i,k}^g v_i^k+\sum_{i\in \mathcal{I}_1} \bar{\lambda}_{i,k}^h \nabla h_i(x^k)
-\sum_{i\in \mathcal{I}_2\cup\mathcal{I}_5^k}\bar{\lambda}_{i,k}^G \nabla G_i(x^k)-\sum_{i\in \mathcal{I}_3\cup\mathcal{I}_6^k}\bar{\lambda}_{i,k}^H \nabla H_i(x^k)+\sum_{i\in L'_k}\alpha_i^k \nu_i^k,\label{AKKT3}
\end{eqnarray}
and the set of vectors $\{v_i^k\}_{i\in \mathcal{I}_4^k}\cup
\{\nabla h_i(x^k)\}_{i\in \mathcal{I}_1}\cup \{\nabla G_i(x^k)\}_{i\in \mathcal{I}_2\cup\mathcal{I}_5^k}\cup \{\nabla H_i(x^k)\}_{i\in \mathcal{I}_3\cup\mathcal{I}_6^k}\cup \{\nu_i^k\}_{i\in L'_k}$ is linearly independent.

Since the index sets are all finite, for every large $k$, we  may assume $\mathcal{I}_4^k\equiv\mathcal{I}_4$, $L'_k\equiv L$, $\mathcal{I}_5^k\equiv\mathcal{I}_5$ and $\mathcal{I}_6^k\equiv\mathcal{I}_6$ without loss of generality. Hence the set of vectors
  $\{v_i^k\}_{i\in \mathcal{I}_4}\cup
\{\nabla h_i(x^k)\}_{i\in \mathcal{I}_1}\cup \{\nabla G_i(x^k)\}_{i\in \mathcal{I}_2\cup\mathcal{I}_5}\cup \{\nabla H_i(x^k)\}_{i\in \mathcal{I}_3\cup\mathcal{I}_6}\cup \{\nu_i^k\}_{i\in L}$ is linearly independent.
From $(\ref{AKKT3})$, the condition $(\ref{sumrule})$ reduces to 
\begin{eqnarray}\label{sumrule1}
&&0 = v_0^k+(x^k-x^*)+\sum_{i\in \mathcal{I}_4} \bar{\lambda}_{i,k}^g v_i^k+\sum_{i\in \mathcal{I}_1} \bar{\lambda}_{i,k}^h \nabla h_i(x^k)
-\sum_{i\in \mathcal{I}_2}\bar{\lambda}_{i,k}^G \nabla G_i(x^k)-\sum_{i\in \mathcal{I}_3}\bar{\lambda}_{i,k}^H \nabla H_i(x^k)\nonumber \\
&&~~~ -\sum_{i\in\mathcal{I}_5} \bar{\lambda}_{i,k}^G \nabla G_i(x^k)
-\sum_{i\in\mathcal{I}_6}\bar{\lambda}_{i,k}^H \nabla H_i(x^k) +\sum_{i\in L}\alpha_i^k \nu_i^k,
\end{eqnarray}
where  $\alpha_i^k> 0$, $i\in L$,  
$\bar{\lambda}_{i,k}^g> 0$ if $i\in \mathcal{I}_4$, $\bar{\lambda}_{i,k}^G>0, \bar{\lambda}_{i,k}^H>0$ or $\bar{\lambda}_{i,k}^G \bar{\lambda}_{i,k}^H=0$ if  $i\in \mathcal{J}^*$.

Step 3:  We now prove that the sequence $\{(
\bar{\lambda}_{k}^g, 
\bar{\lambda}_{k}^h,\bar{\lambda}_{k}^G,\bar{\lambda}_{k}^H, \sum_{i\in L} \alpha_i^k\nu_i^k)\}$ must be bounded. 
To the contrary, assume that it is unbounded. Let $M_k:=\|(
\bar{\lambda}_{k}^g, 
\bar{\lambda}_{k}^h,\bar{\lambda}_{k}^G,\bar{\lambda}_{k}^H)\|+\|\sum_{i\in L} \alpha_i^k\nu_i^k\|$.  Then there exists a subsequence $K$   such that $M_k\to\infty$ and
$$\displaystyle\lim_{k\to\infty,k\in K}\frac{(
\bar{\lambda}_{k}^g, 
\bar{\lambda}_{k}^h,\bar{\lambda}_{k}^G,\bar{\lambda}_{k}^H,\sum_{i\in L} \alpha_i^k\nu_i^k)}{M_k}=(
\lambda^g, 
\lambda^h,\lambda^G, \lambda^H, {\eta^*}),$$ 
where 
{$\eta_i^*=\lim_{k\rightarrow \infty} \frac{\alpha_i^k}{M_k}\eta_i^k$, $i\in L$ and $\eta_i^*=0$, $i\notin L$}.  Since $\eta^k_i \in  \mathcal{N}_{C_i}(x^k_i)$ and $\alpha_i^k/M_k>0$, it follows from the outer semicontinuity of the limiting normal cone that $\eta^*_i\in \mathcal{N}_{C_i}(x^*_i)$  and $\eta^*\in \mathcal{N}_{C}(x^*)$ for each $i\in L$. 
It is easy to see that $\lambda_{i}^g\geq 0$ if $i\in \mathcal{I}_4$ and $\mbox{ either } \lambda_i^G> 0, \lambda_i^H>  0 \mbox{ or }\lambda_i^G \lambda_i^H=0 \ \  \forall i \in \mathcal{J}^*$.
Without loss of generality,  assume that 
$v_i^k\to v_i^*$ and
$v_0^k\to v_0^*$.  Then by the outer semi-continuity of the limiting subdifferential, we have  $ 
v_i^*\in \partial g_i(x^*)$ and $v_0^*\in \partial f(x^*)$.  
Dividing $M_k$ on both sides of \((\ref{sumrule1})\) and taking limits with $k \rightarrow \infty, k\in K$, we have 
\begin{eqnarray*}
0 = \sum_{i\in \mathcal{I}_4} \lambda_i^g v_i^* +\sum_{i\in \mathcal{I}_1} \lambda_i^h \nabla h_i(x^*)-\sum_{i\in\mathcal{I}_2\cup\mathcal{I}_5} \lambda_i^G \nabla G_i(x^*)
-\sum_{i\in\mathcal{I}_3\cup\mathcal{I}_6}\lambda_i^H \nabla H_i(x^*)
+ \eta^*.
\end{eqnarray*}
By  RCPLD(ii), since the vectors  $\lambda^g_{\mathcal{I}_4},\lambda^h_{\mathcal{I}_1},\lambda^G_{\mathcal{I}_2\cup\mathcal{I}_5},\lambda^H_{\mathcal{I}_3\cup\mathcal{I}_6}$, $\eta^*$  are not all equal to zero,  $x^k\rightarrow x^*, v_i^k \rightarrow v_i^*$, $\frac{\alpha^k_i\eta^k_i}{M_k}\rightarrow \eta^*_i, i\in L$, the set of vectors $$\{v_i^k\}_{i\in \mathcal{I}_4}\cup\{\nabla h_i(x^k)\}_{i\in \mathcal{I}_1}\cup
\{\nabla G_i(x^k)\}_{i\in \mathcal{I}_2\cup\mathcal{I}_5}\cup \{\nabla H_i(x^k)\}_{i\in \mathcal{I}_3\cup\mathcal{I}_6}\cup  \{\frac{\alpha^k_i\nu^k_i}{M_k}\}_{i\in L}$$ 
must be  linearly dependent.   Since $\alpha_i^k/M_k>0$ for each $ i\in L$,  this   is a contradiction to the  conclusion in step 2. 
The contradiction proves that $\{(
\bar{\lambda}_{k}^g, 
\bar{\lambda}_{k}^h, \bar{\lambda}_{k}^G,\bar{\lambda}_{k}^H, \sum_{i\in L}\alpha^k_i \nu_i^k)\}$ is bounded.  

Without loss of generality, assume that $(
\bar{\lambda}_{k}^g, 
\bar{\lambda}_{k}^h,\bar{\lambda}_{k}^G,\bar{\lambda}_{k}^H, \sum_{i\in L}\alpha^k_i \nu_i^k)$ $\to (\lambda^g, \lambda^h,\lambda^G, \lambda^H, \sum_{i\in L} \nu_i^*)$ as $k\rightarrow \infty$. It is easy to see that $\mbox{ either } \lambda_i^G> 0, \lambda_i^H>  0 \mbox{ or }\lambda_i^G \lambda_i^H=0 \ \  \forall i \in \mathcal{J}^*$,
 $\displaystyle\lim_{k\to\infty}\displaystyle  \sum_{i\in L}\alpha_i^k\nu_i^k=\sum_{i\in L} \nu_i^*\in \mathcal{N}_{C}(x^*)$. 
Taking limits for a subsequence in $(\ref{sumrule1})$, we derive that $x^*$ is an M-stationary point of $(\rm P)$.
\BOX

We now perform the next task of proving that RCPLD is a sufficient condition for the error bound property under extra regularity conditions. First we prove the following result which shows that RCPLD is persistent locally.

\begin{proposition}\label{persis}  Assume that  $\rm RCPLD$ holds at  a feasible point $x^*$ of the system
\begin{equation}
\begin{array}{l}
g(x)\leq 0, 
 \quad h(x)=0,
\quad x\in C.
 \end{array} \label{F1}
\end{equation} Then  $\rm RCPLD$ holds at every point belongs to a sufficiently small neighborhood  of $x^*$.
\end{proposition}   
{\bf Proof.}  Consider any  sequence $x^k\rightarrow x^*$ as $k\to\infty$.  
It is obvious that
$\rm RCPLD$ (i) holds at each $x^k$ when $k$ is sufficiently large.
Assume that $\{\nabla h_i(x^*)\}_{i\in \mathcal{I}_1}$ with $\mathcal{I}_1\subseteq \{1,\cdots,m\}$  is a basis for span$\{\nabla h_i(x^*)\}_{i=1}^m$.  Then $\{\nabla h_i(x^k)\}_{i\in \mathcal{I}_1}$ is a basis of span$\{\nabla h_i(x^k)\}_{i=1}^m$ for any large $k$.
We now show that  RCPLD (ii) holds at $x^k$ for any sufficiently large $k$. 
To contrary, assume that  RCPLD(ii) does not hold at each point of a subsequence $\{x^k\}_{k\in K_0}$. Then there exist an index set $\mathcal{I}_2^k\subseteq I_g^k:=\{i=1,\cdots,n:g_i(x^k)=0\}$,  $v_i^k\in\partial g_i(x^k)$,  $i\in \mathcal{I}_2^k$ and  a nonzero vector $(\lambda^g_{k},\lambda^h_{k}, \eta^k )$ with $\lambda^g_{\mathcal{I}_2^k,k}\geq 0$, $\eta^k:=(\eta^k_1,\dots, \eta^k_l) \in \mathcal{N}_{C}(x^k)$ such that 
\begin{eqnarray}\label{persisnew}
0 = \sum_{i\in \mathcal{I}_2^k} \lambda_{i,k}^g v_i^k +\sum_{i\in \mathcal{I}_1} \lambda_{i,k}^h \nabla h_i(x^k)+\eta^k,
\end{eqnarray} but the set of vectors
$
 \{v_i^{k,s}\}_{i\in \mathcal{I}_2^k}
 \cup \{\nabla h_i(y^{k,s})\}_{ i\in \mathcal{I}_1}\cup \{\nu_i^{k,s}\}_{i\in L},
 $   where $L:=\{1,\dots, l:\eta_i^{k,s}\neq 0\}$, 
is linearly independent  for some
  sequences
 $y^{k,s}\to x^k$ with $\partial g_i(y^{k,s}) \ni v_i^{k,s}\to v_i^k \  (i\in \mathcal{I}_2^k)$,
$\nu_i^{k,s}:=\{0\}^{s_i}\times \{\eta_i^{k,s}\}\times \{0\}^{t_i}$, $\eta_i^{k,s}\in \mathcal{N}_{C_i}(y_i^{k,s}), \eta_i^{k,s}\to\eta_i^k$ as $s\to\infty$.
Since $\mathcal{I}_2^k\subseteq I_g^k\subseteq I_g^*$ is a finite set, we may consider a subsequence such that $\mathcal{I}_2^k=\mathcal{I}_2$ for every large $k\in K_0$.
{Let $M_k:=\|\lambda^g_k\|+\|\lambda_{k}^h\|+\| \eta^k\|$. Suppose there exists a subsequence $K\subseteq K_0$   such that 
$$\displaystyle\lim_{k\to\infty,k\in K}\frac{(
{\lambda}_{k}^g, 
\lambda_{k}^h,\eta^k)}{M_k}=(
\lambda^g, 
\lambda^h, \eta)$$ with ${\lambda}_i^g\geq 0, i\in \mathcal{I}_2$.
 Without loss of generality,  assume that 
$v_i^k\to v_i\in \partial g_i(x^*)$ and we assume $\displaystyle\lim_{k\to\infty,k\in K}\frac{\eta_i^k}{M_k}=\eta_i\in \mathcal{N}_{C_i}(x_i^*)$, $i=1,\cdots,l$. 
Since $\eta_i^{k,s}\to\eta_i^k$ as $s\to\infty$, 
 we also have $\frac{\eta_i^{k,s}}{M_k}\to \eta_i, i\in L$ as $k,s\to\infty$.}

By the diagonalization law,  there exists a sequence $\{z^k\}$ converging to $x^*$ such that for each $k$, $\bar{v}_i^{k}\in \partial g_i(z^k)$, $\bar{v}_i^{k}\to{v}_i^{k}, i=1,\cdots,n$, $\bar{\eta}^k_i\in \mathcal{N}_{C_i}(z_i^k)$, 
$\frac{\bar{\eta}_i^k}{M_k}\to \eta_i, i\in L$
and $\{\bar{v}_i^{k}\}_{i\in \mathcal{I}_2} \cup \{\nabla h_i(z^k)\}_{ i\in \mathcal{I}_1}\cup\{\bar{\nu}_i^{k}\}_{i\in L}$ is linearly independent for all large $k$ and $\bar{\nu}_i^{k}:=\{0\}^{s_i}\times \{\bar{\eta}_i^k\}\times \{0\}^{t_i}$.

Dividing by $M_k$ in both sides of $(\ref{persisnew})$ 
and letting $k\to \infty, k\in K$, we have that $\eta\in \mathcal{N}_{C}(x^*)$ and
$$
0= \sum_{i\in \mathcal{I}_2} {\lambda}_i^g v_i+\sum_{i\in \mathcal{I}_1} \lambda_i^h \nabla h_i(x^*)+\eta.
$$
From  RCPLD (ii), 
the set of vectors $\{\bar{v}_i^{k}\}_{i\in \mathcal{I}_2} \cup \{\nabla h_i(z^k)\}_{ i\in \mathcal{I}_1}\cup 
\{\frac{\bar{\nu}_i^{k}}{M_k}\}_{i\in L}$  
must be  linearly dependent, which is a  contradiction.
The contradiction shows  that  $\rm RCPLD$ holds at $x^k$ for $k$ sufficiently large.
\BOX


\begin{thm}\label{err1}
Assume the $RCPLD$ holds at a feasible point $x^*$ of the system $(\ref{F1})$,
 $g_i(\cdot), i=1,\cdots,n$ are subdifferentially regular  and $C$ is  Clarke regular around $x^*$, 
then  there exist $\alpha>0$ and $\varepsilon>0$ such that 
\begin{eqnarray*}
d_{{\cal F}_1}(x) \leq \alpha (\|g_+(x)\| + \|h(x)\|),
 \quad \forall x\in \mathbb{B}_{\varepsilon}(x^*)\cap C,
\end{eqnarray*}
where ${\cal F}_1$ denotes the set of feasible points satisfying system (\ref{F1}).
\end{thm}  
{\bf Proof.}  If $x^*$ is an interior point of ${\cal F}_1$, then $d_{{\cal F}_1}(x) =0$ for all $x\in \mathbb{B}_{\varepsilon}(x^*)\cap C$ and hence the result holds automatically. Now assume that $x^*$ lies in the boundary of ${\cal F}_1$.  Let $\phi_1(x):=\|g_+(x)\| + \|h(x)\|$.
We rewrite the feasible set by
$
{\cal F}_1:=\left \{x\in C: \phi_1(x)=0
\right \}.
$ 
Assume for a contradiction that there exists $C \ni x^k\to x^*$ such that
$
d_{{\cal F}_1}(x^k)>k \phi_1(x^k).
$

Obviously $x^k\notin {\cal F}_1$. Let $y^k$ be the projector of $x^k$ to ${\cal F}_1$.  Then $d_{{\cal F}_1}(x^k)=\|y^k-x^k\|\not =0$ and $\displaystyle \lim_{k\to\infty}y^k=x^*$. 
For each $k$, $y^k$ is an optimal solution of the following problem:
\begin{eqnarray*}
({\rm P}'_k)~~~\min && F^k(x):=\|x-x^k\| \\
 {\rm s.t.} && g(x)\leq 0, h(x)=0, x\in C,
\end{eqnarray*}
where $\|\cdot\|$ denotes the 2-norm.
Since the RCPLD persists in a neighborhood of $x^*$,  RCPLD holds at $y^k$ for $k$ sufficiently large.  From Theorem \ref{OC1}, $y^k$ is a limiting stationary point of $({\rm P}'_k)$. By the optimality condition, 
there exist parameters $\lambda_{i,k}^g\geq 0$, $v_i^k\in \partial g_i(y^k)$ for $i\in I(y^k)$ and $\lambda_{i,k}^h, i=1,\cdots,m$ such that 
\begin{eqnarray}\label{kt}
0= \frac{y^k-x^k}{\|y^k-x^k\|}+\sum_{i\in I(y^k)} \lambda_{i,k}^g v_i^k +\sum_{i=1}^m \lambda_{i,k}^h \nabla h_i(y^k)+\eta^k.
\end{eqnarray}
Let $\nu_i^k:=\{0\}^{s_i}\times \{\eta_i^k\}\times \{0\}^{t_i}$, $\sum_{i\in L_k} \nu_i^k=\eta^k \in \mathcal{N}_{C}(y^k)$ and $L_k:=\{1,\cdots,l:\eta_i^k\neq 0\}$.

Assume that $\{\nabla h_i(x^*)\}_{i\in \mathcal{I}_1}$ with $\mathcal{I}_1\subseteq \{1,\cdots,m\}$  is a basis for span$\{\nabla h_i(x^*)\}_{i=1}^m$.
From Lemma \ref{lem3-1}, we obtain
$
\mathcal{I}_2^k\subseteq I(y^k)\cap {\rm supp}(\lambda_{i,k}^g)
$
and 
$L'_k\subseteq L_k$ with 
$\{\bar{\lambda}_{k}^g, \bar{\lambda}_{k}^h, \alpha^k\}$, $\bar{\lambda}_{i,k}^g> 0$ for $i\in \mathcal{I}_2^k$ and $\alpha_i^k> 0$, $i\in L'_k$
  such that 
\begin{eqnarray}\label{AKKT30}
0= \frac{y^k-x^k}{\|y^k-x^k\|}+ \sum_{i\in \mathcal{I}_2^k} \bar{\lambda}_{i,k}^g v_i^k+\sum_{i\in \mathcal{I}_1} \bar{\lambda}_{i,k}^h \nabla h_i(y^k)+
\sum_{i\in L'_k}\alpha_i^k \nu_i^k
\end{eqnarray}
with bounded multipliers $\{(\bar{\lambda}_{k}^g, \bar{\lambda}_{k}^h, \xi^k)\}$ from Theorem \ref{OC1}, for $k\to\infty$, $\xi^k:=\sum_{i\in L'_k}\alpha_i^k \nu_i^k\in \mathcal{N}_{C}(y^k)$.
Let $M\geq\|(\bar{\lambda}_{k}^g, \bar{\lambda}_{k}^h)\|_1$ for sufficiently $k$.

From the subdifferentially regularity of $g_i(\cdot), i=1,\cdots,n$, for sufficiently large $k$, 
\begin{eqnarray*}
g_i(x^k) - g_i(y^k) -\langle v_i^k, x^k-y^k \rangle +\frac{1}{4M}\|x^k-y^k\|\geq 0.
\end{eqnarray*}
Similarly, $\langle \xi^k,x^k-y^k \rangle \leq \frac{1}{4}\|y^k-x^k\|$.
Furthermore, for each $i=1,\cdots,m$,
\begin{eqnarray*}
h_i(x^k) =h_i(y^k) +\langle \nabla h_i(y^k), x^k-y^k \rangle +o(\|x^k-y^k\|).
\end{eqnarray*}

Then from \((\ref{AKKT30})\), we have 
\begin{eqnarray*}
 \|y^k-x^k\|&=&\langle \sum_{i\in \mathcal{I}_2^k} \bar{\lambda}_{i,k}^g v_i^k+\sum_{i\in \mathcal{I}_1} \bar{\lambda}_{i,k}^h \nabla h_i(y^k)+\xi^k,x^k-y^k \rangle\\
&\leq & \sum_{i\in \mathcal{I}_2^k} \bar{\lambda}_{i,k}^g (g_i(x^k)-g_i(y^k)) +\sum_{i\in \mathcal{I}_1} \bar{\lambda}_{i,k}^h (h_i(x^k)-h_i(y^k)) \\
&& +\frac{1}{4}\|y^k-x^k\| + (\sum_{i\in \mathcal{I}_2^k} \bar{\lambda}_{i,k}^g +\sum_{i\in \mathcal{I}_1} |\bar{\lambda}_{i,k}^h|)\frac{1}{4M}\|y^k-x^k\|\\
&\leq &  \sum_{i\in \mathcal{I}_2^k} \bar{\lambda}_{i,k}^g g_i(x^k)+\sum_{i\in \mathcal{I}_1} \bar{\lambda}_{i,k}^h h_i(x^k)
+\frac{1}{2}\|y^k-x^k\|.
\end{eqnarray*}
This means
\begin{eqnarray*}
d_{{\cal F}_1}(x^k)= \|y^k-x^k\|&\leq& 2M (\sum_{i\in \mathcal{I}_2^k} \max\{0, g_i(x^k)\}+\sum_{i\in \mathcal{I}_1} |h_i(x^k)|)\\
&\leq&2M(\|\max\{0,g(x^k)\}\|_1+\|h(x^k)\|_1),
\end{eqnarray*}
which is a contradiction.
Therefore the error bound property holds.
\BOX

The error bound property for the general system (\ref{feasibilitys}) can be now obtained  from Theorem \ref{err1}. It extends \cite[Theorem 5.1]{gzl} to allow nonsmooth inequality constraints and abstract constraints.
\begin{corollary}\label{gerr1}
Assume the $RCPLD$ holds at $x^*\in \mathcal{F}$,
 $g_i(\cdot), i=1,\cdots,n$ are subdifferentially regular and $C$ is Clarke regular around $x^*$, the constraint  $(G(x) ,H(x)) \in \Omega^p$ satisfies the strict complementarity condition at $x^*$,
then  there exist $\alpha>0$ and $\varepsilon>0$ such that 
\begin{eqnarray}
d_{{\cal F}}(x) \leq \alpha \phi(x),
 \quad \forall x\in \mathbb{B}_{\varepsilon}(x^*)\cap C, \label{Errorb}
\end{eqnarray}
where  $\phi(x):=\|g_+(x)\| + \|h(x)\| + \sum_{i\in {\cal I}^*} |G_i(x)|+ \sum_{i\in {\cal K}^*} |H_i(x)|$ and ${\cal F}$ denotes the set of feasible points satisfying system (\ref{feasibilitys}).
\end{corollary} 
{\bf Proof.} Since the strict complementarity holds  at $x^*$, we have $\mathcal{J}^*=\emptyset$. Hence for all $x\in {\cal F}$ sufficiently close to $x^*$, we can represent it as a solution to the system
$$g_i(x)\leq 0,\ i=1,\cdots,n,\\
 h_i(x)=0,\ i=1,\cdots,m, G_i(x)=0, i\in \mathcal{I}^*, H_i(x)=0 , i\in  \mathcal{K}^*, x\in C.$$ 
 Then the error bound property follows by Theorem \ref{err1} and the equivalence of the finite dimensional norm.
\BOX

\section{Sufficient conditions for RCPLD}
Note that although the RCPLD is a weak condition, it may not be easy to verify. In this section we investigate sufficient conditions for RCPLD which are stronger but easier to verify.

It is easy to see that the following  well-known constraint qualification implies RCPLD. 
\begin{defn}\label{nna*} Let $x^*\in {\cal F}$. We say that the no nonzero abnormal multiplier constraint qualification (NNAMCQ) holds at $x^*$ if there is no nonzero vector
$(\lambda^g,\lambda^h, \lambda^G,\lambda^H, \eta^*) \in \mathbb{R}^{n}\times\mathbb{R}^{m}\times \mathbb{R}^{p}\times\mathbb{R}^p\times\mathbb{R}^d$ satisfying $\lambda_i^g\geq 0$ for $i\in I_g^*$, 
and $\mbox{either } \lambda_i^G> 0, \lambda_i^H>  0 \mbox{ or }\lambda_i^G \lambda_i^H=0, \forall i \in \mathcal{J}^*$,  $\eta^*\in \mathcal{N}_{C}(x^*)$, $v^*_i \in \partial g_i(x^*)$ for $i\in I_g^*$  such that 
\begin{eqnarray*}
&&0 = \sum_{i\in I_g^*} \lambda_i^g v_i^* +\sum_{i=1}^m \lambda_i^h \nabla h_i(x^*)
 -\sum_{i\in\mathcal{I}^*\cup\mathcal{J}^*} \lambda_i^G \nabla G_i(x^*)
 -\sum_{i\in\mathcal{K}^*\cup\mathcal{J}^*}\lambda_i^H \nabla H_i(x^*)
 + \eta^*.\nonumber
\end{eqnarray*} 
\end{defn}

 It is  {obvious} that when the rank of $\{\nabla h_i(x^*)\}_{i=1}^m\cup \{\nabla G_i(x^*)\}_{i\in \mathcal{I}^*}\cup \{\nabla H_i(x^*)\}_{i\in \mathcal{K}^*}$  is equal to $d$,   RCPLD holds automatically. This condition is easy to verify and  moreover it is not just a constraint qualification but also a sufficient condition for error bounds without imposing any regularity conditions. 
\begin{thm}\label{err2}For a feasible point $x^*\in \mathcal{F}$,
suppose that 
the rank of $ \{\{\nabla h_i(x^*)\}_{i=1}^m\cup \{\nabla G_i(x^*)\}_{i\in \mathcal{I}^*}\cup \{\nabla H_i(x^*)\}_{i\in \mathcal{K}^*}\}$  is equal to $d$. Then the error bound property (\ref{errorbp}) holds at $x^*$. 
\end{thm}  
{\bf Proof.}  
 Around the point $x^*$, we can equivalently formulate the complementarity system $(G(x) ,H(x)) \in \Omega^p$ as 
$$G_i(x)=0, i\in \mathcal{I}^*, H_i(x)=0 , i\in  \mathcal{K}^*, (G_i(x),H_i(x))\in \Omega, i\in \mathcal{J}^*.$$
 Hence the constraints
$G_i(x)=0, i\in \mathcal{I}^*, H_i(x)=0 , i\in  \mathcal{K}^*$ can be treated as equality constraints. We denote by $\phi(x):=\max\{0,g_i(x), i=1,\cdots,n, 
 |h_i(x)|,i=1,\cdots,m,\ |G_i(x)|, i\in \mathcal{I}^*, |H_i(x)|,  i\in  \mathcal{K}^*,\ d_{\Omega}(G_i(x),H_i(x)), i\in \mathcal{J}^*\}$.
 
 If $x^*$ is an interior point of ${\cal F}$, then $d_{{\cal F}}(x) =0$ for all $x\in \mathbb{B}_{\varepsilon}(x^*)\cap C$ and hence the result holds automatically. Now assume that $x^*$ lies in the boundary of ${\cal F}$.  
We rewrite the feasible set by
$
{{\cal F}:=\left \{x\in C: \phi(x)=0
\right \}.}
$ 
To a contrary, assume that there exists $C \ni x^k\to x^*$ such that
\begin{eqnarray}\label{contrerr}
d_{{\cal F}}(x^k)>k \phi(x^k).
\end{eqnarray}

For any $i=1,\cdots,m$, $h_i(x^*)=0$ and by  the Taylor expansion,
\begin{equation}
h_i(x^k) =h_i(x^*) +\langle \nabla h_i(x^*), x^k-x^* \rangle +o(\|x^k-x^*\|). 
\label{Taylor}
\end{equation}
From $(\ref{contrerr})$, $\|x^*-x^k\|\geq d_{{\cal F}}(x^k)>k \phi(x^k)\geq k |h_i(x^k)|$, which implies that $\displaystyle \lim_{k\to\infty} \frac{h_i(x^k)}{\|x^*-x^k\|}=0$. Taking subsequence if necessary, let $d^*:=\displaystyle \lim_{k\to\infty}\frac{x^k-x^*}{\|x^k-x^*\|}$. Dividing the both sides of (\ref{Taylor}) by $\|x^k-x^*\|$ and letting $k\to\infty$, we have
$
\langle \nabla h_i(x^*), d^* \rangle=0.
$
Similarly from the above discussion, we have
$
\langle \nabla G_i(x^*), d^* \rangle=0,\ \forall i\in \mathcal{I}^*,
\quad \langle \nabla H_i(x^*), d^* \rangle=0,\ \forall i\in \mathcal{K}^*.
$
This means that $d^*$ is linearly independent with all vectors in $\mbox{ span } \{\{\nabla h_i(x^*)\}_{i=1}^m\cup \{\nabla G_i(x^*)\}_{i\in \mathcal{I}^*}\cup \{\nabla H_i(x^*)\}_{i\in \mathcal{K}^*}\}$. Since by assumption  the rank of $\mbox{span } \{\{\nabla h_i(x^*)\}_{i=1}^m\cup \{\nabla G_i(x^*)\}_{i\in \mathcal{I}^*}\cup \{\nabla H_i(x^*)\}_{i\in \mathcal{K}^*}\}$ is $d$, this  is impossible since these vectors lie in $\mathbb{R}^d$. The proof of the theorem is therefore complete.
\BOX

In the following definition we extend the well-known concept of  RCRCQ {(see \cite{M-S} for the smooth equality and inequality systems and \cite[Definition 3.4]{gly1} for system with complementarity constraints)} to our general system (\ref{feasibilitys}). RCRCQ is  stronger than the RCPLD but may be easier to verify.
\begin{defn} Let $x^*\in {\cal F}$. 
We say that  the relaxed constant rank constraint qualification $(\rm RCRCQ)$ holds  at $x^*$ if for all sufficiently large $k$, 
 any  index sets $\mathcal{I}_4\subseteq I_g^*$, $\mathcal{I}_5\subseteq \mathcal{J}^*$,  $\mathcal{I}_6\subseteq  \mathcal{J}^*$, {$\mathcal{L}\subseteq \{1,\cdots,l\}$} and  any vectors $v^*_i\in \partial g_i(x^*)(i\in \mathcal{I}_4)$, $ \eta^*_i\in \mathcal{N}_{C_i}(x^*_i)$ ($i\in \mathcal{L}$), $ \nu_i^*:=\{0\}^{s_i}\times \{\eta_i^*\}\times \{0\}^{t_i}$ with $s_i:=q_1+\cdots +q_{i-1}$ and $t_i:=q_{i+1}+\cdots + q_l$,  the set of vectors 
$$\{v_i^*\}_{i\in \mathcal{I}_4}\cup \{\nabla h_i(x^*)\}_{i=1}^m\cup \{\nabla G_i(x^*)\}_{i\in \mathcal{I}^*\cup \mathcal{I}_5}\cup\{\nabla H_i(x^*)\}_{i\in \mathcal{K}^*\cup\mathcal{I}_6}\cup 
 \{ \nu_i^*\}_{i\in \mathcal{L}},
$$ and the set of vectors 
$$\{v_i^k\}_{i\in \mathcal{I}_4}\cup \{\nabla h_i(x^k)\}_{i=1}^m\cup \{\nabla G_i(x^k)\}_{i\in \mathcal{I}^*\cup \mathcal{I}_5}\cup\{\nabla H_i(x^k)\}_{i\in \mathcal{K}^*\cup\mathcal{I}_6}\cup 
 \{ \nu_i^k\}_{i\in \mathcal{L}},
$$
where  $x^k\not =x^*$, 
{$\nu_i^k:=\{0\}^{s_i}\times \{\eta_i^k\}\times \{0\}^{t_i}$,}
 have the same rank for all  sequences $\{x^k\}, \{v^k\}, \{\eta^k\}$ satisfying $x^k\rightarrow x^*$, $v^k\rightarrow v^*$, $\eta^k\rightarrow \eta^*$ as $k\rightarrow \infty$,  $v_i^k \in \partial g_i(x^k)$,  $\eta_i^k \in \mathcal{N}_{C_i}(x_i^k)$.
\end{defn}

\begin{proposition}  RCRCQ  implies RCPLD.
\end{proposition}
{\bf Proof.} Let $x^*\in {\cal F}$ satisfy RCRCQ. Taking $\mathcal{I}_4\cup \mathcal{I}_5\cup \mathcal{I}_6\cup \mathcal{L}=\emptyset$, we have  that RCPLD (i) holds. We now show the RCPLD (ii) holds at $x^*$. Let  $\mathcal{I}_1\subseteq \{1,\cdots,m\}$, $\mathcal{I}_2\subseteq \mathcal{I}^*$,  $\mathcal{I}_3\subseteq \mathcal{K}^*$ be such that  the set of vectors $\{\nabla h_i(x^*)\}_{i\in \mathcal{I}_1}\cup \{\nabla G_i(x^*)\}_{i\in \mathcal{I}_2}\cup \{\nabla H_i(x^*)\}_{i\in \mathcal{I}_3}$ 
 is a basis for  
 $ \mbox{ span }\large \{ \{\nabla h_i(x^*)\}_{i=1}^m\cup \{\nabla G_i(x^*)\}_{i\in \mathcal{I}^*}\cup \{\nabla H_i(x^*)\}_{i\in \mathcal{K}^*} \large \}.$

 Let  $\mathcal{I}_4\subseteq I_g^*$, $\mathcal{I}_5,\mathcal{I}_6\subseteq \mathcal{J}^*$ be the index sets such that 
 there exists a nonzero vector $(\lambda^g,\lambda^h, \lambda^G,\lambda^H, \eta^*)$ satisfying $\lambda_i^g\geq 0$ for $i\in \mathcal{I}_4$, 
and $\mbox{either } \lambda_i^G> 0, \lambda_i^H>  0 \mbox{ or }\lambda_i^G \lambda_i^H=0, \forall i \in \mathcal{J}^*$,  ${\eta^*=(\eta_1^*,\cdots,\eta_l^*)}\in \mathcal{N}_{C}(x^*)$, $v^*_i \in \partial g_i(x^*)$ for $i\in {\cal I}_4$  satisfying
\begin{eqnarray*}
&&0 = \sum_{i\in \mathcal{I}_4} \lambda_i^g v_i^* +\sum_{i\in \mathcal{I}_1} \lambda_i^h \nabla h_i(x^*)
 -\sum_{i\in\mathcal{I}_2\cup\mathcal{I}_5} \lambda_i^G \nabla G_i(x^*)
 -\sum_{i\in\mathcal{I}_3\cup\mathcal{I}_6}\lambda_i^H \nabla H_i(x^*)
 + \eta^*.\nonumber
\end{eqnarray*} 
Then the vectors 
\begin{eqnarray*}
&&\{v_i^* \}_{i\in \mathcal{I}_4} \cup \{\nabla h_i(x^*)\}_{i\in \mathcal{I}_1}  
 \cup \{\nabla G_i(x^*)\}_{i\in\mathcal{I}_2\cup\mathcal{I}_5} \cup \{\nabla H_i(x^*)\}_{i\in\mathcal{I}_3\cup\mathcal{I}_6}\cup  \{ \nu_i^*\}_{i\in L},
\end{eqnarray*} 
 where $ \nu_i^*:=\{0\}^{s_i}\times \{\eta_i^*\}\times \{0\}^{t_i}$  and $L:=\{1,\dots, l:\eta_i^*\neq 0\}$, 
is linearly dependent.  Since $\{\nabla h_i(x^*)\}_{i\in \mathcal{I}_1}\cup \{\nabla G_i(x^*)\}_{i\in \mathcal{I}_2}\cup \{\nabla H_i(x^*)\}_{i\in \mathcal{I}_3}$ 
 is a basis for  
 $ \mbox{ span }\large \{ \{\nabla h_i(x^*)\}_{i=1}^m\cup \{\nabla G_i(x^*)\}_{i\in \mathcal{I}^*}\cup \{\nabla H_i(x^*)\}_{i\in \mathcal{K}^*} \large \},$ 
 it follows that \begin{eqnarray*}\label{rcr}
&&\{v_i^* \}_{i\in \mathcal{I}_4} \cup \{\nabla h_i(x^*)\}_{i=1} ^m
 \cup \{\nabla G_i(x^*)\}_{i\in\mathcal{I}^*\cup\mathcal{I}_5} \cup \{\nabla H_i(x^*)\}_{i\in\mathcal{K}^*\cup\mathcal{I}_6}\cup  \{ \nu_i^*\}_{i\in L} 
\end{eqnarray*} 
is linearly dependent as well. 
By RCRCQ, it follows that 
\begin{eqnarray*}
\{v_i^k \}_{i\in \mathcal{I}_4} \cup \{\nabla h_i(x^k)\}_{i=1} ^m 
 \cup \{\nabla G_i(x^k)\}_{i\in\mathcal{I}^*\cup\mathcal{I}_5} \cup \{\nabla H_i(x^k)\}_{i\in\mathcal{K}^*\cup\mathcal{I}_6}\cup \{ \nu_i^k\}_{i\in L} 
 \end{eqnarray*} 
is linearly dependent 
 for all  sequences $\{x^k\}, \{v^k\}, \{\eta^k\}$ satisfying $x^k\rightarrow x^*$, $v^k\rightarrow v^*$, $\eta^k\rightarrow \eta^*$,  $v_i^k \in \partial g_i(x^k)$,  $\eta_i^k \in \mathcal{N}_{C_i}(x_i^k)$, $\nu_i^k:=\{0\}^{s_i}\times \{\eta_i^k\}\times \{0\}^{t_i}$. 
From RCPLD (i), $\{\nabla h_i(x^k)\}_{i\in \mathcal{I}_1}\cup \{\nabla G_i(x^k)\}_{i\in \mathcal{I}_2}\cup \{\nabla H_i(x^k)\}_{i\in \mathcal{I}_3}$ is a basis of $\mbox{ span }\{ \{\nabla h_i(x^k)\}_{i=1}^m\cup \{\nabla G_i(x^k)\}_{i\in \mathcal{I}^*}\cup \{\nabla H_i(x^k)\}_{i\in \mathcal{K}^*}\}$ and thus
 \begin{eqnarray*}
\{v_i^k \}_{i\in \mathcal{I}_4} \cup 
\{\nabla h_i(x^k)\}_{i\in\mathcal{I}_1}  
 \cup \{\nabla G_i(x^k)\}_{i\in\mathcal{I}_2\cup\mathcal{I}_5} \cup \{\nabla H_i(x^k)\}_{i\in\mathcal{I}_3\cup\mathcal{I}_6}\cup \{ \nu_i^k\}_{i\in L}
 \end{eqnarray*} 
 is linearly dependent. Since  $L\subseteq L_k:=\{1,\dots, l:\eta_i^k\neq 0\}$ for any sufficiently large $k$, RCPLD (ii) holds.
\BOX

\begin{defn}  We say that the Linear Constraint Qualification (LCQ) holds if
all functions $g_i,h_i, G_i,H_i$ are linear and the set $C$ is the union of  finitely many polyhedral sets. 
\end{defn}
 We now show that LCQ implies RCRCQ.
\begin{proposition} LCQ implies RCRCQ holds at each $x^*\in {\cal F}$.
\end{proposition}
{\bf Proof.}  
Since $C$ is the union of finitely many polyhedral convex sets, for all $x^k \rightarrow x^*$ and sufficiently large $k$, we have
 $\eta^k=\eta^*$ if $\eta^k\rightarrow \eta^*$ and $\eta^k\in \mathcal{N}_{C}(x^k)$. Moreover all functions are linear. Therefore the matrix $J'(x,v_{{\cal I}_2}, \nu_L):=\left [\begin{array}{ccccc}
 v_{{\cal I}_2} &   \nabla h(x)^T&   \nabla G_{{\cal I}^*\cup {\cal I}_3}(x)^T
& \nabla H_{{\cal K}^*\cup {\cal I}_4}(x)^T &
\nu_L \end{array}
\right ] $ is a constant matrix for all   $x, v_{{\cal I}_2}=\nabla g_{{\cal I}_2}(x), \eta\in \{\eta^k, \eta^*\}$, $\nu_i:=\{0\}^{s_i}\times \{\eta_i\}\times \{0\}^{t_i}$, $i\in L$ and therefore RCRCQ holds.
\BOX

\begin{eg}\label{nonsmoothex} Consider the nonsmooth system:
\begin{eqnarray*}
&& h(x):=2x_1+x_2=0,\\
&& g(x):= x_1 +x_2-\max\{\frac{1}{2}, x_3\}-x_4+1\leq 0,\\
&& x\in C:=C_1\times C_2,
\end{eqnarray*}
\end{eg}
where $C_1$ is the graph of a  continuous function $\varphi:[-1,1]\to\mathbb{R}$ as shown in Fig.1  and 
$C_2:=\{(t,1-t): t\in [0,1]\}$. We set $\varphi(x)=0$ when $x=0, 2^{-n}, -2^{-n}$, for $n=0,1,2,\cdots$.
Between any two points of the form $2^{-n-1}$ and $2^{-n}$, or $-2^{-n}$ and $-2^{-n-1}$, the graph of $\varphi$ describes the edge of an isosceles triangle whose
apex is located at $(2^{-n-2}+2^{-n-1},1)$ or $(-2^{-n-2}-2^{-n-1},1)$.
Consider the feasible solution $x^*=(x^*_1,x^*_2,x^*_3,x^*_4)=(0,0,\frac{1}{2},\frac{1}{2})$.
\begin{figure}
  \centering
    \includegraphics[width=0.6\textwidth]{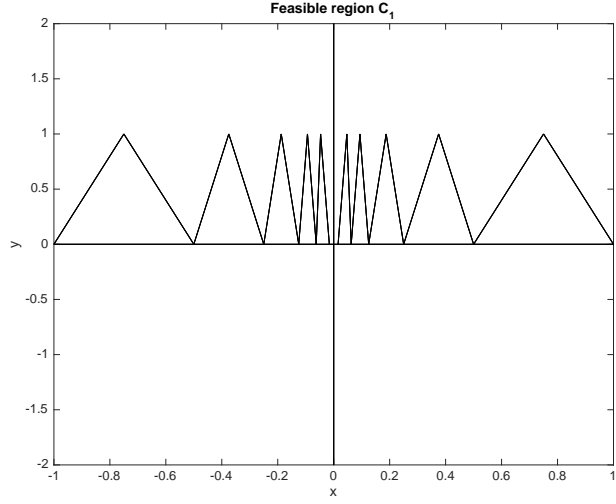}
    \caption{Feasible region $C_1$}
\end{figure}
We have $\nabla h(x)=(2,1,0,0) $, 
$$\partial g(x)=\left\{\begin{array}{ll}
(1,1,-1,-1)  & {\rm if}\ x_3>\frac{1}{2},\\ 
\{(1,1,-1,-1) ,(1,1,0,-1)   \} & {\rm if}\ x_3=\frac{1}{2},\\ 
(1,1,0,-1)  & {\rm if}\ x_3<\frac{1}{2}.
\end{array}\right.
$$
Since $C_1$ is  {symmetrical}, we only give the expression of the normal cone for the case when $(x_1,x_2)\in C_1$ and  $x_1\geq 0$, for $n=0,1,2,\cdots$,
\begin{eqnarray*}
\lefteqn{\mathcal{N}_{C_1}(x_1,x_2)=}\\
&&\left\{
(\eta_1,\eta_2) :
\begin{array}{ll}
\eta_1=0, \eta_2\in \mathbb{R}& {\rm if}\ x_1\in (2^{-n-1},2^{-n}), x_2=0,\\
\eta_2=-2^{-n-2}\eta_1, \eta_1\in \mathbb{R} & {\rm if}\ x_1\in (2^{-n-1},2^{-n-2}+2^{-n-1}), x_2>0,\\
\eta_2=2^{-n-2} \eta_1, \eta_1\in \mathbb{R}& {\rm if}\ x_1\in (2^{-n-2}+2^{-n-1},2^{-n}), x_2>0,\\
\eta_2=\alpha \eta_1,  \alpha=\pm 2^{-n-2}\  {\rm or}\ \eta_2>0\ {\rm with}\\
\alpha  \in (-\infty,-2^{-n-2})\cup(2^{-n-2},+\infty)
 & {\rm if}\ x_1=2^{-n-2}+2^{-n-1}, x_2=1,\\
\eta_1=0, {\rm or}\
\eta_2=\alpha \eta_1, \alpha=-2^{-n-2}, 2^{-n-3} 
& {\rm if}\ x_1=2^{-n-1}, x_2=0,\\
\eta_1 \eta_2=0 & {\rm if}\ x_1=0, x_2=0,\\
\end{array}\right.
\end{eqnarray*}
and
\begin{eqnarray*}
\mathcal{N}_{C_2}(x_3,x_4)=\left\{(\eta_3,\eta_4) :
\begin{array}{ll}
 \eta_4\geq \eta_3,& {\rm if}\ x_3=0, x_4=1\\
\eta_3= \eta_4,& {\rm if}\ x_3 \in (0,1), x_4=1-x_3\\
\eta_3\geq \eta_4,& {\rm if}\ x_3=1, x_4=0
\end{array}\right \}.
\end{eqnarray*}
Also we have $\mathcal{N}_{C}(x)= \mathcal{N}_{C_1}(x_1,x_2)\times\mathcal{N}_{C_2}(x_3,x_4)$. 
Hence $\mathcal{N}_{C}(x^*)= \{(\eta_1,\eta_2):\eta_1\eta_2=0\}\times \{(\alpha,\alpha):\alpha\in\mathbb{R}\}$.
Let $\eta^*=(0,0, 0,0)$ being  an element in the normal cone $\mathcal{N}_{C}(x^*)$.  For any $x^k\to x^*$ with $x_1^k=x_2^k=0, x_3^k>\frac{1}{2}$, $v^k=(1,1,-1,-1)\in \partial g(x^k)$, $ \eta^k=(0,0,\alpha_k,\alpha_k)\in {\cal N}_C(x^k)$ with $\alpha_k \rightarrow 0$,  we have
$$
 J'(x^k,v^k,  \eta^k):=\left [\begin{array}{ccc}
  v^k &  \nabla h(x^k) &     \eta^k \end{array}
\right ]=
\left [\begin{array}{ccc}
1 &2 &0 \\
1 & 1 &0 \\
-1 & 0 & \alpha_k\\
-1 &0&\alpha_k
\end{array}\right ].$$ 
$r(J'(x^k,v^k,  \eta^k))=3$ if $\alpha_k\neq 0$ but $r(J'(x^*,v^*,  \eta^*))=2$. This shows that RCRCQ  fails at $x^*$.

Since the condition 
\begin{eqnarray}\label{ex4-1-1}
0= \lambda v^*+\mu\nabla h({x^*}) +\eta^*
\end{eqnarray}
holds when 
$\eta^* = (1,0, 1,1) \in \mathcal{N}_{C}(x^*)$, $v^*=(1,1,-1,-1)  \in \partial g(x^*)$, $\lambda=1=-\mu$. Hence the NNAMCQ fails at $x^*$. 

We now verify that RCPLD holds. Actually it is easy to see that  $$0=\lambda v^*+ \mu\nabla h(x^*) + \eta^*, \quad v^*\in \partial g(x^*), \quad \eta^*\in \mathcal{N}_{C}(x^*) $$ holds with  $\mu, \lambda\geq 0, \eta^*$ not all equal to zero if and only if $v^*=(1,1,-1,-1)$ and 
either (a) or (b) below holds:
\begin{itemize}
\item[{\rm (a)}]
$\lambda=-\mu$,
${\eta}^*=\mu (-1,0,-1,-1) $.
\item[{\rm (b)}]
$\lambda=-2\mu$,
${\eta}^*=2\mu(0,-\frac{1}{2}, 1,1) $.
\end{itemize}
In either case (a) or (b), the set $L=\{1,2\}$. For any $x^k\to x^*=(0,0,\frac{1}{2},\frac{1}{2})$, 
denote by $\eta^k:=(\eta_1^k,\eta_2^k,\alpha_k,\alpha_k) $ where $(\eta_1^k,\eta_2^k)\in \mathcal{N}_{C_1}(x_1^k,x_2^k)$, $(\alpha_k,\alpha_k)\in \mathcal{N}_{C_2}(x_3^k,x_4^k)$ and $\eta^k \rightarrow  \eta^*$.  We also denote by $\nu_1^k:=(\eta_1^k,\eta_2^k,0,0), \nu_2^k:=(0,0, \alpha_k,\alpha_k)$. 
Note that at any $x^k$ such that  $x^k\rightarrow x^*, v^k\rightarrow v^*, v^k\in \partial g(x^k)$ and $k$ is large enough, $g(x)$ is differentiable and $v^k=\nabla g(x^k)=(1,1,-1,-1)$. 
Since  the matrix
$$J'(x^k, v^k,\nu_L^k)=\left [\begin{array}{cccc}
1 &2 & \eta_1^k&0 \\
1 & 1& \eta_2^k &0 \\
-1 & 0 & 0 & \alpha_k\\
-1 &0& 0&\alpha_k
\end{array}\right ]$$
 is not full rank,
the set of vectors $\nabla h(x^k),v^k$, $\nu_1^k$, $\nu_2^k$  is always linearly dependent and thus RCPLD holds at $x^*$.

Note that 
if $\alpha_k-\eta_1^k+2\eta_2^k\neq 0$ for large $k$,
it is easy to see that the rank of the matrix
$$\left [\begin{array}{ccc}
1 &2 & \eta_1^k \\
1 & 1&\eta_2^k  \\
-1 & 0  & \alpha_k\\
-1 &0&\alpha_k
\end{array}\right ]$$
equals to $3$ and thus
the set of vectors $\nabla h(x^k), v^k$, $\eta^k$ is linearly independent. Hence  as pointed out in Remark \ref{remark}, our definition for RCPLD is weaker than the condition that the set of vectors $\nabla h(x^k), v^k $, $\eta^k$ is linearly dependent.

\section{Applications to bilevel programs}
In this section we apply RCPLD to the combined program (CP).
Throughout this section we assume that the  value function $V(x)$ is  Lipschitz continuous at the  point of interest. For the case where $Y(x)=Y$ is independent of $x$, we can use Danskin's theorem and for the general case one can use  Proposition \ref{partialV}  which is a special case of  \cite[Theorem 6.5.2]{c}.
Other weaker sufficient conditions for Lipschitz continuity of the value function as well as the estimates for its subdifferentials can be found for examples in \cite{glyz}.
\begin{proposition}[Danskin's Theorem](\cite[page 99]{clsw} or \cite{Danskin}) \label{Danskin}  Let $Y\subseteq \mathbb{R}^s$ be a compact set and $f(x,y)$ be a function defined on $\mathbb{R}^d\times \mathbb{R}^s$ that is continuously differentiable at $x^*$. Then the value function 
$V(x):=\min \{ f(x,y): y\in Y\}$ is Lipschitz continuous near $x^*$ and its Clarke subdifferential is 
$\partial^c V(x^*) = co \{ \nabla_x f(x^*, y): y\in S(x^*)\}$,
where $\partial^c V(x^*)=co \partial V(x^*) $ is the Clarke subdifferential of $V$ at $x^*$. 
\end{proposition}
\begin{proposition}\label{partialV}Assume that the set-valued map $Y(x)$ is uniformly bounded around $x^*$, i.e., there exists a neighborhood $\mathbb{U}$ of $x^*$ such that the set $\bigcup_{x\in \mathbb{U}} Y(x)$ is bounded. Suppose that MFCQ holds at $y$ for all $y\in S(x^*)$. Then the valued function $V(x)$ is Lipschitz continuous near $x^*$ and 
$$\partial^c V(x^*) \subseteq co W(x^*),$$
where 
\begin{eqnarray*}
W(x^* )&:=& \bigcup_{y\in S(x^*) } \left \{\nabla_x f(x^*, y)+u \nabla_x g(x^*,y)+v \nabla_x h(x^*,y): (u,v)\in M(x^*, y) \right \},\\
M(x^*,y)&:=& \left\{ (u,v): 
\begin{array}{c}
0=\nabla_y f(x^*,y)+u\nabla_y g(x^*,y)+v \nabla_y h(x^*,y)\\
u\geq 0,\quad  \langle g(x,y),u\rangle=0
\end{array}
\right\},
\end{eqnarray*}
where $u\nabla_y g:=\sum_{i=1}^m u_i \nabla_y g_i$.
\end{proposition}
Note that if in addition to the assumptions of Proposition \ref{partialV}, {$S$ is inner semicontinuous at $(x^*,y^*)$ for some $y^*\in S(x^*)$,} then the union $\bigcup_{y\in S(x^*) }$ sign can be omitted in Proposition \ref{partialV}; see \cite[Corollary 1.109]{m1}.
 In the  case  where the LICQ holds at each $y\in S(x^*)$, the set of multipliers $M(x^*, y)$ is a singleton and by  
\cite[Corollary 5.4]{GauvinD}, the inclusion becomes an equality in the above proposition 
and $-V(x)$ is Clarke regular around $x^*$. {In this case, LICQ for the lower level problem holds at every $y, y\in S(x)$, for all $x$ near $ x^*$  due to the outer semi-continuity of the solution mapping $S(x)$. Thus we have
$$\partial_x (f-V)(x,y)=  \partial_x^c (f-V)(x,y)=\nabla_
x f(x,y)- \partial^c V(x) = \nabla_x f(x,y)- co W(x).$$
By Carath\'{e}odory's theorem, the convex set $co W(x)\subseteq \mathbb{R}^d$ can be represented by not more than $d+1$ elements. It follows that  for any $w\in co W(x)$, there exist $\mu^i\geq 0, \sum_{i=1}^{d+1}\mu^i=1$, $M(x, y^i)=\{(u^i,v^i)\}, y^i\in S(x)$ such that
$w=\sum_{i=1}^{d+1}\mu^i(\nabla_x f(x, y^i)+u^i \nabla_x g(x,y^i)+v^i \nabla_x h(x,y^i) )$.}

Given a feasible vector $(x^*,y^*,u^*,v^*)$ of the problem (CP),  define the following index sets:
 \begin{eqnarray*}
 && I_G^*=I_G(x^*,y^*):=\{i: G_i(x^*,y^*)=0\},\\
 &&{\cal I}^*={\cal I}(x^*,y^*,u^*):=\{i: g_i(x^*,y^*)=0, u_i^* >0\},\\
 && {\cal J}^*={\cal J}(x^*,y^*,u^*):=\{i: g_i(x^*,y^*)=0, u_i^* =0\},\\
 && {\cal K}^*={\cal K}(x^*,y^*,u^*):=\{i: g_i(x^*,y^*)<0, u_i^*=0\}.
\end{eqnarray*}

Theorem \ref{OC1} can now be applied to the problem (CP) to obtain the M-stationary condition under RCPLD at any local optimal solution.  
\begin{thm} \label{rcpb} Let $(x^*,y^*,u^*,v^*)$ be a local solution of $(\rm CP)$ and suppose that the value function $V(x)$ is Lipschitz continuous at $x^*$.
If the RCPLD holds at $(x^*,y^*,u^*,v^*)$, then
$(x^*,y^*,u^*,v^*)$ is an M-stationary point of problem (CP). 
\end{thm}

In the rest of this section, we apply the sufficient conditions for RCPLD which are introduced in  Section 4 to the bilevel programs. The following theorem follows immediately from Theorem \ref{err2}. The condition is easy to verify since the nonsmooth constraint $f(x,y)-V(x)\leq 0$ is not needed in the verification.

\begin{thm} \label{thm5.2} Let $(x^*,y^*,u^*,v^*)$ be a local solution of $(\rm CP)$ and suppose that the value function $V(x)$ is Lipschitz continuous at $x^*$.
If the rank of the matrix 
$$J^*=\left [\begin{array}{cccc}
\nabla  (\nabla_y f+u\nabla_y g+v\nabla_y h)(x^*,y^*) & \nabla_y h(x^*,y^*)^T & \nabla_{y} g_{{\cal I}^*\cup{\cal J}^*}(x^*,y^*)^T \\
\nabla h(x^*,y^*)  &0 &0\\
\nabla H(x^*,y^*)  &0&0 \\
\nabla g_{{\cal I}^*}(x^*,y^*)  & 0&0
\end{array}\right ]$$
is equal to $d+s+m+n-|{\cal K}^*|$,
then 
{RCPLD holds and}
$(x^*,y^*,u^*,v^*)$ is an M-stationary point of problem (CP). Moreover $(x^*,y^*,u^*,v^*)$ is a local optimal solution of the penalized problem for some $\mu\geq 0$:
\begin{eqnarray*}
({\rm CP}_\mu)~~~~~~\min_{x,y,u,v} && F(x,y) +\mu \phi_{CP}(x,y,u,v)
\end{eqnarray*}
where $\phi_{CP}(x,y,u,v):=(f(x,y)-V(x))_++\|H(x,y)\|+\|h(x,y)\|+\|G_+(x,y)\|+\|\nabla_y L (x,y,u,v)\|+\sum_{i=1}^md_{\Omega}(-g_i(x,y),u_i))$.
\end{thm}
{\bf Proof.} Let $
\Gamma^*:=\left[\begin{array}{c}
\nabla_y g_{{\cal K}^*}(x^*,y^*)^T\\
0
\end{array}\right]$.
It is easy to see that 
$$
r\left(\left[\begin{array}{cc}
J^*& \Gamma^*\\
0& I_{|\mathcal{K}^*|}
\end{array}\right]\right)
=r(J^*)+|\mathcal{K}^*|=d+s+m+n,
$$ 
where $I_{|\mathcal{K}^*|}$ is the identity matrix of size $|\mathcal{K}^*|$.
From Theorem \ref{err2}, the error bound property holds for the problem (CP), i.e., there exist $\alpha>0$ and $\varepsilon>0$ such that
\begin{eqnarray*}
d_{{\cal F}_{CP}}(x,y,u,v) \leq \alpha\phi_{CP}(x,y,u,v), \quad \forall (x,y,u,v)\in \mathbb{B}_{\varepsilon}(x^*,y^*,u^*,v^*),
\end{eqnarray*}
where ${\cal F}_{CP}$ denotes the feasible region of problem (CP).
It follows from Clarke's exact penalty principle \cite[Proposition 2.4.3]{c} that 
the problem $({\rm CP}_\mu)$ is exact with $\mu\geq L_F\alpha$ with $L_F$ being the Lipschitz constant of the function $F$.
\BOX

\begin{corollary}\label{corcpcp}
Let $(x^*,y^*,u^*,v^*)$ be a local solution of $(\rm CP)$. Suppose that the value function $V(x)$ is Lipschitz continuous at $x^*$ and LICQ holds at $y^*$.
Suppose 
the  matrix
$$SJ^*=\left [\begin{array}{l}
\nabla h(x^*,y^*)\\ 
\nabla H(x^*,y^*)\\ 
\nabla g_{{\cal I}^*}(x^*,y^*)\end{array}\right ]
 $$
has full column rank  $d+s$.   Then $(x^*,y^*,u^*,v^*)$ is an M-stationary point of problem (CP) and a local optimal solution of the penalized problem $({\rm CP}_\mu)$ for some $\mu\geq 0$.
\end{corollary}  
{\bf Proof.} Since LICQ holds at $y^*$, the rank of the matrix 
$$SJ_1^*=\left [\begin{array}{cc}
\nabla_y h(x^*,y^*)^T& \nabla_{y} g_{{\cal I}^*\cup{\cal J}^*}(x^*,y^*)^T 
\end{array}\right ]$$
equals to $|{\cal I}^*|+|{\cal J}^*|+n$. Then $r(J^*)=r(SJ^*)+r(SJ_1^*)=d+s+m+n-|{\cal K}^*|$ and thus the conclusions in Theorem \ref{thm5.2} hold.
\BOX

The following example illustrates the application of Theorem \ref{thm5.2}  and Corollary \ref{corcpcp}.
\begin{eg}\label{exyz} Consider the following bilevel program:
\begin{eqnarray*}
\min && F(x,y)\\
{\rm s.t.} && x\in [-3,2],\\
&& H(x,y):=x^2+y-2=0,\\
&& y\in S(x):=\argmin\limits_{y}\ f(x,y):=y^3-3y\\
&&~~~~~~~~~~~~~~~~~~~~{\rm s.t.}\  g_1(x,y):=x-y\leq 0,\\
&&~~~~~~~~~~~~~~~~~~~~~~~~~g_2(x,y):=y-3\leq 0,
\end{eqnarray*}
where $F(x,y)$ is a Lipschitz continuous function.
It is easy to see that the solution set for the lower level program is
\begin{eqnarray*}
S(x)= \left \{ \begin{array}{ll}
\{x\} & \mbox{ if } x\in[-3,-2)\cup (1,2],\\
\{-2,1\} &\mbox{ if } x=-2,\\
\{1\} & \mbox{ if } x\in(-2,1],
\end{array} \right.
\end{eqnarray*}
and the value function is
\begin{eqnarray*}
V(x)= \left \{ \begin{array}{ll}
x^3-3x & \mbox{ if } x\in[-3,-2)\cup (1,2],\\
-2 &\mbox{ if }  x\in[-2,1].
\end{array} \right.
\end{eqnarray*}
In fact since {the lower level feasible set $Y(x):=\{y\in \mathbb{R}:x\leq y\leq 3\}$ is uniformly bounded whenever $x\in [-3,2]$ and LICQ holds at each $y\in S(x)$, $x\in [-3,2]$,}  we can conclude that the value function is Lipschitz continuous without actually calculating it.
For any $x$, the KKT condition for the lower level problem is
$$0=3y^2-3-u_1+u_2,\  g_i(x,y)\leq 0,\ u_i\geq 0,\ u_i g_i(x,y)=0, i=1,2.$$ 
Hence the combined problem can be written as follows:
\begin{eqnarray*}
\min_{  x,y, u}  && F(x,y)\\
{\rm s.t.} && f(x,y)-V(x)\leq 0,\\
&& H(x,y)=0,\\
&& \nabla_y L(x,y,u):=3y^2-3-u_1+u_2=0,\\
&& x\in [-3,2],  (-g(x,y),u)\in \Omega^2.
\end{eqnarray*}
It is easy to see that the feasible region of the bilevel problem contains three points: $(x^*,y^*)=(-2,-2), (\hat{x},\hat{y})=(1,1), (\tilde{x},\tilde{y})=(-1,1)$.  From calculation,
\begin{eqnarray*}
\nabla (\nabla_y L)( x, y, u)=(0,6y,-1,1), \nabla H(x,y)=(2x,1), \nabla g_1(x,y)=(1,-1), \nabla g_2
(x,y)=(0,1).
 \end{eqnarray*}

 Suppose that the optimal solution of the bilevel problem is $(x^*,y^*)=(-2,-2)$ and the one for the combined program is $(x^*,y^*,u_1^*,u_2^*)=(-2,-2,9,0)$. 
The index sets ${\cal I}^*=\{1\}, {\cal K}^*=\{2\}$ and ${\cal J}^*$ is empty. The rank of the matrix
 $$SJ^*=\left [\begin{array}{c}
\nabla H(x^*,y^*)  \\
\nabla g_1(x^*,y^*)  
\end{array}\right ]
=
\left [\begin{array}{cc}
-4 & 1 \\
1 & -1  
\end{array}\right ]$$
is equal to $2$. Hence by Corollary \ref{corcpcp}, $(x^*,y^*,u_1^*,u_2^*)$ is an M-stationary point of problem (CP). 

{ Suppose that 
the optimal solution of the bilevel problem is $(\tilde{x},\tilde{y})=(-1,1)$ and the one for the combined program is $(\tilde{x},\tilde{y},\tilde{u}_1,\tilde{u}_2)=(-1,1,0,0)$. 
The index sets $\tilde{\cal K} =\{1,2\}$ and both $\tilde{\cal I}, \tilde{\cal J}$ are empty. The rank of the matrix
 $$J(\tilde{x},\tilde{y},\tilde{u}_1,\tilde{u}_2)=\left [\begin{array}{c}
\nabla (\nabla_y f+u\nabla_y g)(\tilde{x},\tilde{y}) \\
\nabla H(\tilde{x},\tilde{y}) 
\end{array}\right ]
=
\left [\begin{array}{cc}
0 &6  \\
-2 & 1 
\end{array}\right ]$$
is equal to $2$. Hence by Theorem \ref{thm5.2}, $(\tilde{x},\tilde{y},\tilde{u}_1,\tilde{u}_2)$ is an M-stationary point of problem (CP).}
Moreover 
$(x^*,y^*,u_1^*,u_2^*)$ and $(\tilde{x},\tilde{y},\tilde{u}_1,\tilde{u}_2)$ are local solutions of the penalized problem for some $\mu\geq 0$:
\begin{eqnarray*}
\min_{  x,y, u}  &&F(x,y)+\mu \left ((f(x,y)-V(x))_++ |H(x,y)|+|\nabla_y L(x,y,u)|+|\min\{-g(x,y),u\}|\right )\\
{\rm s.t.}  && x\in [-3,2].
\end{eqnarray*}

\end{eg}

Since RCRCQ is a stronger condition for RCPLD, the following result follows from 
Theorem \ref{OC1} and Propositions \ref{Danskin} and \ref{partialV}.
\begin{thm}\label{thm5.1} 
Let $(x^*,y^*,u^*,v^*)$ be a local solution of $(\rm CP)$. Suppose that either $Y(x)=Y$ is independent of $x$ with  $Y$ compact or the set-valued map $Y(x)$ is uniformly bounded around  $x^*$ and MFCQ holds at $y$ for all $y\in S(x^*)$. 
Given 
 index sets ${\cal I}_3, {\cal I}_4\subseteq {\cal J}^*$  and ${\cal I}_2 \subseteq  I_G^*$,  $\alpha\in \{0,1\}$, denote  the matrix
\begin{eqnarray*}
 J'(x,y,u,v,\alpha, w):=
 \left [\begin{array}{ccc}
\nabla(\nabla_{y} f+u\nabla_y g+v\nabla_y h)(x,y) &\nabla_y h(x,y)^T&\nabla_y g(x,y)^T \\
  \nabla h(x,y)  & 0& 0 \\
  \nabla H(x,y) &0& 0\\
   \nabla G_{{\cal I}_2}(x,y)&0&0\\
   \nabla g_{{\cal I}^*\cup {\cal I}_3}(x,y)&0&0 \\
 0 & 0& E_{\mathcal{K}^*\cup{\cal I}_4}  \\
    \alpha ( \nabla f(x,y)- (w,0) )&0&0 \end{array}\right ],
\end{eqnarray*}
 where $E_{\mathcal{K}^*\cup{\cal I}_4}\subseteq \mathbb{R}^{(|\mathcal{K}^*|+|{\cal I}_4|)\times m}$ denotes the matrix with $e_i^T$ as its rows, $i\in \mathcal{K}^*\cup{\cal I}_4$ and $e_i\in \mathbb{R}^m$ is the vector such that the $i$-th component is one and others are zero. 
Assume that  the matrix $J'(x^*,y^*,u^*,v^*,\alpha, w^*)$  where $w^* \in coW(x^*)$ has the same rank with the matrix  $J'(x^k,y^k,u^k,v^k,\alpha, w^k)$  for all sequences $\{x^k\}, \{y^k\}, \{u^k\}, \{v^k\}, \{w^k\}$ satisfying $x^k\rightarrow x^*, y^k\rightarrow y^*, u^k\rightarrow u^*, v^k\rightarrow v^*, w^k\rightarrow w^*$,  $w^k\in coW(x^k)$.
Then $\rm RCRCQ$ for $(CP)$ hold at $(x^*,y^*,u^*,v^*)$ and $(x^*,y^*,u^*,v^*)$ is an M-stationary point of problem (CP). 
\end{thm}

The following example illustrate the result.
\begin{eg}\label{exp} Consider the following bilevel program:
\begin{eqnarray*}
\min && F(x,y)\\ 
{\rm s.t.} && H(x,y):=x_1-x_2+y-\frac{1}{2}=0,\\
&&y\in \argmin_y\ f(x,y):=x_1 \exp(y)-x_2 \exp(y)\\
&&~~~~~~~~~~{\rm s.t.}\ g_1(y):=-y-\ln2\leq 0,\\
&&~~~~~~~~~~~~~~~g_2(y):=y-\ln2\leq 0,\,
\end{eqnarray*}
where $F(x,y)$ is a Lipschitz continuous function.
It is easy to see that the solution set for the lower level program is
\begin{eqnarray*}
S(x)= \left \{ \begin{array}{ll}
[-\ln2,\ln2] &\mbox{ if } x_1=x_2,\\
\{\ln2\}& \mbox{ if } x_1<x_2,\\
\{-\ln2\} &\mbox{ if } x_1>x_2,
\end{array} \right.
\end{eqnarray*}
and the value function 
\begin{eqnarray*}
V(x)= \left \{ \begin{array}{ll}
0 & \mbox{ if } x_1=x_2,\\
2(x_1-x_2)&\mbox{ if }  x_1<x_2,\\
\frac{1}{2}(x_1-x_2) & \mbox{ if } x_1>x_2.
\end{array} \right.
\end{eqnarray*}
In fact we do not need the above explicit representation of the value function. Indeed,
since the constraint set $Y=[-\ln 2, \ln2]$ is compact, by Danskin's theorem,  for any $x$ around $x^*$ the Clarke subdifferential of the value function is equal to 
$$
\partial^c V(x)=coW(x)= co\{\nabla_x f(x,y_x)| y_x\in S(x)\}=\displaystyle co\{\exp(y_x)| y_x\in S(x)\}(1,-1).
$$
The combined problem can be written as follows:
\begin{eqnarray*}
\min_{  x,y, u}  && F(x,y)\\
{\rm s.t.} && f(x,y)-V(x)\leq 0,  H(x,y)=0,\\
&& \nabla_y L(x,y,u):=x_1 \exp(y)-x_2 \exp(y)-u_1+u_2=0,\\
&&  (-g(y),u)\in \Omega^2.
\end{eqnarray*}
It is easy to see that the bilevel program only has three kinds of optimal solutions, $(x_1^*,x_2^*,y^*)=(a,a,\frac{1}{2})$ with $a\in \mathbb{R}$, $(\hat{x}_1,\hat{x}_2,\ln2)$ with $\hat{x}_1 <\hat{x}_2$ and $(\tilde{x}_1,\tilde{x}_2,-\ln2)$  with $\hat{x}_1 >\hat{x}_2$ .  
We now verify that  RCRCQ holds in the following three cases.

(a) Consider the optimal solution $(x_1^*,x_2^*,y^*)=(a,a,\frac{1}{2})$ with $a\in \mathbb{R}$, the corresponding solution for the combined program is $(x^*,y^*,u_1^*,u_2^*)=(a,a,\frac{1}{2},0,0)$. 
Obviously, 
the index sets ${\cal K}^*=\{1,2\}$ and both ${\cal I}^*, {\cal J}^*=\emptyset$. 
For $\alpha\in \{0,1\}$, denote by
 $$ J'(x,y,u,\alpha, \omega):=\left [\begin{array}{ccccc}
\exp(y)&-\exp(y) & x_1 \exp(y)-x_2 \exp(y) & -1 &1\\
1&-1& 1&0&0\\
0& 0& 0&1&0\\
0&  0& 0&0&1\\
{\alpha ( \exp(y) -w)}&{-\alpha ( \exp(y) -w)} &\alpha( x_1 \exp(y)-x_2 \exp(y))&0&0
\end{array}\right ],
$$
where $w\in co\{\exp(y_x)| y_x\in S(x)\}$.  It is easy to see that $r(J'(x^*,y^*,u^*,\alpha, w^*))=r(J'(x,y,u,\alpha, w))=4$ for any $\alpha\in \{0,1\}$,  and any $(x,y,u)$ and  $w$. Hence
RCRCQ holds by Theorem \ref{thm5.1}.

(b) Consider the optimal solution $(\hat{x}_1,\hat{x}_2,\ln2)$ with $\hat{x}_1 <\hat{x}_2$, the corresponding solution for the combined program is $(\hat{x}_1,\hat{x}_2,\ln 2,0,\hat{u}_2)$ with $\hat{x}_1 <\hat{x}_2$ and $\hat{u}_2>0$. 
The index sets $\hat{\cal K}=\{1\}$ and $\hat{\cal I}=\{2\}, \hat{\cal J}=\emptyset$,  and for any $x$ around $\hat{x}=(\hat{x}_1,\hat{x}_2)$, $S(x)=\{\ln2\}$,  $coW(x)=(2,-2)^T$.
For $\alpha\in \{0,1\}$, denote by
  $$ J'(x,y,u,\alpha, w):=\left [\begin{array}{ccccc}
\exp(y)&-\exp(y) & x_1 \exp(y)-x_2 \exp(y) & -1 &1\\
1&-1& 1&0&0\\
0& 0& 1&0&0\\
0&  0& 0&1&0\\
\alpha(\exp(y)-2)&-\alpha(\exp(y)-2)&\alpha( x_1 \exp(y)-x_2 \exp(y))&0&0
\end{array}\right ].
$$
It is easy to see that $r(J'(x,y,u,\alpha, w))=4$ or any $\alpha\in \{0,1\}$ and any $(x,y,u)$ and $w$. Hence
RCRCQ holds from Theorem \ref{thm5.1}.

(c) Consider the optimal solution $(\tilde{x}_1,\tilde{x}_2,-\ln2)$ with $\tilde{x}_1 >\tilde{x}_2$, the corresponding solution for the combined program is $(\tilde{x}_1,\tilde{x}_2,-\ln 2,\tilde{u}_1,0)$ with $\tilde{x}_1 >\tilde{x}_2$ and $\hat{u}_1>0$. Similarly as in case (b), RCRCQ holds.

\end{eg}

\section*{Acknowlegement} The authors would like to thank  the anonymous referees for their helpful suggestions and comments which help us to improve the presentation of the paper.

\end{document}